\documentclass[draft]{amsart}

\usepackage{graphics}
\usepackage[all]{xy}

\usepackage{amssymb,amsmath}
\usepackage{psfrag}

\usepackage{mathrsfs} 
\usepackage{amsfonts}
\usepackage{amscd}

\usepackage{graphicx}
\usepackage{epstopdf}
\usepackage{tikz}
\usetikzlibrary{arrows,automata}

\usepackage{graphicx,epsfig,float}

\usepackage{url}

\makeatletter
\def\url@leostyle{%
  \@ifundefined{selectfont}{\def\UrlFont{\sf}}{\def\UrlFont{\small\ttfamily}}}
\makeatother
\newtheorem{thm}{Theorem}[section]
\newtheorem{lem}[thm]{Lemma}
\newtheorem{prop}[thm]{Proposition}
\newtheorem{cor}[thm]{Corollary}
\newtheorem{fact}[thm]{Fact}

\newtheorem{clm}[thm]{Claim}

\newcommand{\pf}{{\bf Proof. }}

\renewcommand{\bar}{\overline}

\newcommand{\bJ}{{\mathbf J}}

\newcommand{\cl}{\mathrm{cl}}
\newcommand{\id}{\mathrm{id}}

\newcommand{\dcl}{\mathrm{dcl}}
\newcommand{\dcleq}{\mathrm{dcl}^{\textrm{eq}}}

\begin{document}

\title{On definable Skolem functions and trichotomy}

\author{Bruno Dinis}

\address{ Departamento de Matem\'atica\\
Escola de Ci\^encias e Tecnologia,  Universidade de \'Evora\\
Rua Rom\~ao Ramalho, 59\\ 
7000--671 \'Evora, Portugal}

\email{bruno.dinis@uevora.pt}

\author {M\'ario J. Edmundo}

\address{ Departamento de Matem\'atica\\
Faculdade de Ci\^encias,  Universidade de Lisboa\\
Campo Grande, Edif\'icio C6\\ 
P-1749-016 Lisboa, Portugal}
\email{mjedmundo@fc.ul.pt}

\date{\today}
\thanks{Both authors were supported by FCT - Funda\c{c}\~ao para a Ci\^encia e a Tecnologia, under the projects UIDB/04561/2020 and UIDP/04561/2020 and   the research center  CMAFcIO - Centro de Matem\'atica, Aplica\c{c}\~oes Fundamentais e Investiga\c{c}\~ao Operacional. The first author was also supported by CIMA - Centro de Investigação em Matemática e Aplicações -- UIDP/04674/2020.  \newline
 }


\begin{abstract}
In this paper we give an explicit characterization of o-minimal structures with definable Skolem functions/definable choice. Such structures are, after naming finitely many elements from the prime model, a union of finitely many trivial points each defined over $\emptyset $ and finitely many open intervals each a union of a $\emptyset $-definable family of group-intervals with fixed positive elements.
\end{abstract}

\maketitle

\begin{section}{Introduction}\label{section intro}

In several papers topological invariants of definable sets, such as the fundamental group and  cohomology, were developed in the setting of o-minimal structures ${\mathbb M}=(M,<, (c)_{c\in {\mathcal C}}, $ $ (f)_{f\in  {\mathcal F}}, (R)_{R\in {\mathcal R}})$ with definable Skolem functions/definable choice ( \cite{DEM}, \cite{Edal17}, \cite{ep1}, \cite{ep2}, \cite{ep3}), extending similar work in o-minimal expansions of ordered fields (\cite{bao2}, \cite{bao1}, \cite{bo1}, \cite{bo2}, \cite{dk4}, \cite{dk5},  \cite{eo}, \cite{ew}) and in o-minimal expansions of ordered groups (\cite{EdPa}, \cite{eep}, \cite{ejp1}, \cite{el}, \cite{es}). A referee of a previous version of \cite{DEM}, challenged us to properly discuss the assumption of definable Skolem functions/definable choice and how it differs from the setting of products of finitely many definable group-intervals  which is used in \cite{Edal17}. In fact, the referee asks the following questions: (i) are all examples of such structures, which are not obtained from finitely many group-intervals, similar to the expansion of the real field by the function $f(x,y)=x+y$ with $|y|<1$? (ii) is there a definable family of local group intervals covering all points of $M$? (iii) are all but finitely many points in $\mathbb{M}$ non-trivial?  

Taking into account the questions above we tried to understand the relationship between definable Skolem functions/definable choice and the trichotomy theorem (\cite{pst}). The outcome of this research is the following:

\begin{thm}\label{thm main}
Let $\mathbb{M}$ be an o-minimal structure and let $\mathbb{P}\preceq \mathbb{M}$ be the prime model. Then the following are equivalent:
\begin{enumerate}
\item[(1)]
$\mathbb{M}$ has  definable Skolem functions. 
\item[(2)]
After naming finitely many elements of $\mathbb{P}$, there is a finite collection $\mathcal{Y}$ of $\emptyset $-definable open subintervals of $M$ such that:
\begin{itemize}
\item[-]
$M\setminus \bigcup \mathcal{Y}$ is a finite set of trivial points each defined over $\emptyset $.
\item[-]
each $Y\in \mathcal{Y}$ is the union of a uniformly $\emptyset $-definable family of group-intervals, each with a fixed positive element, parametrized by the end points of the intervals.
\item[-]
each $Y\in \mathcal{Y}$ has  a fixed $\emptyset $-definable element.
\end{itemize}
\item[(3)] $\mathbb{M}$  has definable choice.
\end{enumerate}
\end{thm}

The proof of our main theorem relies on geometrical arguments and constructions. This justifies the development of some technical tools (Section \ref{section auxiliary lemmas}) and the careful analysis, in  Section \ref{section technical lemmas}, of certain definable families of curves obtained from the existence of definable Skolem functions.

Our main theorem gives a full characterization of o-minimal structures with definable Skolem functions/definable choice. This is particularly relevant because so far only loose examples were known such as o-minimal expansions of ordered groups (\cite{vdd}), structures obtained from a product of finitely many definable group-intervals  (\cite{pst}) or the example mentioned in question (i) above.
One can ask for a similar characterization  of o-minimal structures with (parametric) elimination of imaginaries. We return to this question again in Section \ref{section rmks}.

The structure of the paper is the following. For the readers convenience,  we include, in the first two sections,  basic definitions and properties. Most of these are well known or can be found in the literature but are included for completeness. In Section \ref{section prelimp def skolem} we recall the definition of definable Skolem functions/definable choice as well as some properties of these notions. In Section \ref{section prelim trichotomy} we recall the trichotomy theorem (\cite{pst}), the definition of definable group-intervals (\cite{epr}) as well as  some properties of trivial/non-trivial o-minimal structures.  

In Section \ref{section auxiliary lemmas} we prove some auxiliary lemmas. Section \ref{section technical lemmas} is devoted to proving the main technical lemma (Lemma \ref{lem main 0})  required to obtain the main theorem (Theorem \ref{thm main 0}), proved in Section \ref{section main results}. Some final remarks, including an unparametric version of the main theorem, are left to Section \ref{section rmks}.

\medskip
\noindent
{\bf Acknowledgements:} We thank the anonymous referee of a  previous version of \cite{DEM} who's questions challenged us to do this paper!

\end{section}

\begin{section}{Preliminaries on definable Skolem functions}\label{section prelimp def skolem}
\noindent
Recall that an o-minimal structure ${\mathbb M}=(M,<,\cdots)$ has {\it definable Skolem functions}  if and only if for every uniformly definable family $\{X_t\}_{t\in T}$  of nonempty definable subsets of some $M^k$, there is a definable function $h:T\to M^k$ such that:
\begin{itemize}
\item[-]
$h(t)\in X_t$ for all $t\in T$.
\end{itemize}
One says that ${\mathbb M}$ has {\it definable choice} if in addition the definable function $h$ can be chosen so that:
\begin{itemize}
\item[-]
for all $t,t'\in T,$ if $X_t=X_{t'},$ then $h(t)=h(t')$.\\
\end{itemize}

Obviously definable choice implies the existence of definable Skolem functions. In this section we collect a couple of other well known observations about these two notions. 

\begin{lem}\label{lem def skole/choice 0-def and dim one}
For ${\mathbb M}$ to have definable Skolem functions (resp. definable choice) it is enough to consider $\emptyset $-definable families $\{Y_s\}_{s\in S}$ of nonempty definable subsets of $M$.
\end{lem}

\pf
Assume that ${\mathbb M}$ has definable Skolem functions (resp. definable choice) for $\emptyset $-definable families of nonempty definable subsets of $M^k$. Let $\{X_t\}_{t\in T}$ be a definable family of nonempty definable subsets of $M^k$ defined over $b\in M^l$. By varying the parameter $b$ we get  a $\emptyset $-definable family $\{Y_{(t,z)}\}_{(t,z)\in S}$ of nonempty subsets of $M^k$ such that $T=\{t: (t,b)\in S\}$ and $Y_{(t,b)}=X_t$ for all $t\in T$. 
So  there is a definable function $g:S\to M^k$ such that for all $(t,z)$ we have $g(t,z)\in Y_{(t,z)}$.  Therefore, $h:T\to M^k$ given by $h(t)=g(t,b)$ is a definable function such that for all $t\in T$ we have $h(t)\in X_t=Y_{(t,b)}$. Furthermore, if $g$ satisfies the  definable choice property, then $h$ will also satisfy it.

Now assume that ${\mathbb M}$ has definable Skolem functions (resp. definable choice) for $\emptyset $-definable families of nonempty definable subsets of $M$. We prove by induction on $k$ that the same holds for $\emptyset $-definable families of nonempty definable subsets of $M^k$. Let $\{Y_s\}_{s\in S}$ be such a $\emptyset$-definable family. 
By the induction hypothesis, there is $g:T\to M^{k-1}$ a definable Skolem function for $\{\pi (Y_s)\}_{s \in S}$ where $\pi :M^k\to M^{k-1}$ is the projection onto the first $k-1$ coordinates. On the other hand, by the first part of the proof, there is $f:S\to M$ a definable Skolem function for the definable family $\{Z_s\}_{s\in S}$ where $Z_s=\{z\in M: (g(s),z)\in Y_s\}$. It follows that  $h:S\to M^k$ given by $h(s)=(g(s), f(s))$ is a definable Skolem function for $\{Y_s\}_{s\in S}.$ Moreover, if $g$ and $f$ satisfy the  definable choice property, then $h$ will also satisfy it. \qed\\

\begin{prop}\label{prop def skolem/choice elementary substruture}
If ${\mathbb  M}$ has definable Skolem functions (resp. definable choice), then every ${\mathbb K}\prec {\mathbb M}$ has definable Skolem functions (resp. definable choice). 
\end{prop}

\pf
Let ${\mathbb K}\prec {\mathbb M}$. 
Let $\{Y_s\}_{s\in S}$ be a uniformly $\emptyset$-definable family of nonempty definable subsets of $K.$ Then there is a definable function $h:S(M)\to M$ such that for all $z\in S(M)$ we have $h(z)\in Y_z.$ Let $b$ be the parameter in $M^q$ over which $h$ is defined. Then there is a $\emptyset $-definable function $g:S\times D\to K$  such that $b\in D(M)$ and $h(z)=g(z,b)$ for all $z\in S(M).$  Let 
$$D'=\{d\in D:\,\,\,g(s,d)\in Y_s \,\,\,\textrm{for all $s\in S$}\}.$$  
Then $D'$ is $\emptyset $-definable and, since $b\in D'(M)$ and $\mathbb{K}\prec \mathbb{M}$, there is  $c\in D'.$ It follows that if  $f:S\to K$ is given by $f(s)=g(s,c)$ for all $s\in S,$ then $f$ is ${\mathbb K}$-definable and $f(s)\in Y_s$ for all $s\in S.$  

If $h$ satisfies the  definable choice property, then we can choose $g$ satisfying that property and so $f$ will also satisfy it. \qed \\

The following is similar:

\begin{prop}\label{prop def skolem/choice elementary expansion}
If ${\mathbb  M}$ has definable Skolem functions (resp. definable choice), then every ${\mathbb M}\prec {\mathbb K}$ has definable Skolem functions (resp. definable choice). 
\end{prop}

\medskip

Recall that for $C\subseteq M$, the {\it definable closure over $C$} is defined as
$${\rm dcl}(C)=\{u\in M: \{u\}\,\,\,\textrm{is definable over $C$}\}.$$
Equivalently, $u\in \dcl(C)$ if and only if there is a $\emptyset $-definable function $f:Z\subseteq M^l\to M$ and $c\in C^l\cap Z$ such that $f(c)=u.$ Recall also that, since we have an order in $M$, definable closure over $C$ coincides with the {\it algebraic closure over $C$}, defined as 
\[
u\in {\rm acl}(C)\,\,\ \textrm{if and only if $u$ belongs to a finite set definable over $C$}.
\]
As usual if $a=(a_1, \ldots, a_l)\in M^l$ and $C\subseteq M,$ then $aC$ denotes $\{a_1,\ldots, a_l\}\cup C.$\\

From now on we let  ${\mathbb P}\prec {\mathbb M}$ be the prime model  over the empty set  and we let $P$ be its domain. See \cite[Theorem 5.1]{PiSt86} for its existence in the o-minimal context. \\

The following is implicit in the observations before the proof of \cite[Chapter 6, (1.2)]{vdd}:

\begin{prop}\label{prop def skolem dcl}
Assume that ${\mathbb M}$ is sufficiently saturated. Then ${\mathbb M}$ has definable Skolem functions if and only if every nonempty definable subset $X\subseteq M$ defined with parameters $a\in M^l$ has an element in $\dcl(a\,\mathbb{P})$.
\end{prop}

\pf
($\Rightarrow $) Let $X$ be a nonempty definable subset of $M$ defined with parameters $a\in M^l.$ Then there is a uniformly $\emptyset $-definable family $\{X_t\}_{t\in T}$ of nonempty definable subsets of  $P$ such that $a\in T(M)$ and $X_{a}=X.$ By Proposition \ref{prop def skolem/choice elementary substruture}, there is  a ${\mathbb P}$-definable Skolem function  $f:T\to P$ for the family $\{X_t\}_{t\in T}$ of the form $f=g(-,c)$ for some $\emptyset $-definable function $g:T\times D\to P$ and some $c\in P^q$. Then $f(a)\in X\cap {\rm dcl}(a{\mathbb P})$.   

($\Leftarrow $)  Let $\{X_t\}_{t\in T}$ be a uniformly $\emptyset $-definable family of nonempty definable subsets of $M$. Suppose that $T\subseteq M^l$ and fix $v\in T$. Since $X_v$ is defined over $v$, there is an element in   $X_v\cap {\rm dcl}(v\,{\mathbb P})$. So there is a $\emptyset$-definable function $g_v:Z_v\subseteq T\times M^{q_v}\to M$   and $c_v\in P^{q_v}$ such that $(v,c_v)\in Z_v$ and $g_v(v,c_v)\in X_v$.  Let $B_v=\{t\in T: (t,c_v)\in Z_v\}$. Then each $B_v$ is defined over ${\mathbb P}$.

Since $T=\bigcup _v B_v$, by saturation, there are finitely many $v_0,\ldots , v_m$ such that $T=B_{v_0}\cup \ldots \cup B_{v_m}$. For $i=0,\ldots, m$, let  $f_i: B_{v_i}\to M$ be given by $f_i(t)=g_{v_i}(t,c_{v_i})$. Then, by recursion on $m$, there is a ${\mathbb P}$-definable function $f:T\to M$ such that $f(t)\in X_t$ for all $t\in T$. 
\qed \\

\begin{cor}\label{cor prime and dcl}
${\mathbb M}$ has definable Skolem functions if and only if every  $C\subseteq M$, we have ${\rm dcl}(C\, \mathbb{P})\preceq \mathbb{M}$. 
\end{cor}

\pf
It is clear that ${\rm dcl}(C\, \mathbb{P})$ is a substructure of $\mathbb{M}$. To conclude we may assume that ${\mathbb M}$ is sufficiently saturated and use Proposition \ref{prop def skolem dcl} and Tarski's lemma for elementary substructure. 
\qed \\

There is an analogue of Proposition \ref{prop def skolem dcl} for definable choice and definable closure in ${\mathbb M}^{{\rm eq}}$. To this end, let us recall a couple of definitions from \cite[Section 7]{epr}. Given definable sets $X$ and $Y$ and a definable equivalence relation $E$  on $X,$  one says  that a function $f:X/E\to Y$ is {\it $A$-definable} if and only if the set $\{(x,y)\in X\times Y:f([x])=y\}$ is $A$-definable. Equivalently,  $f:X/E\to Y$ is $A$-definable if and only if there is an $A$-definable function $g:X\to Y$ such that   $f([x])=g(x)$ for all $x\in X$ i.e. making the following diagram commutative: 
\begin{equation*} 
\xymatrix{
X \ar[d]_{} \ar[r]^g & Y \\
X/E\ar[ur]^{f} & 
}
\end{equation*}

One says that a subset $U$ of $X/E$ is {\it definable} if it is the image of a definable subset of $X$ under the quotient map $X\to X/E.$
 
 With these definitions, for $a\in X/E$ and $C\subseteq M$ one  sets
$$\dcl (aC)=M\cap \dcleq (aC).$$
Equivalently, $u\in \dcl(aC)$ if and only if there is a $\emptyset $-definable function $f:Z\subseteq X/E\times M^q\to M$ and $c\in C^q\subseteq M^q$ such that $f(a,c)=u.$\\

If $X$ is a nonempty definable subset of $M$ defined with parameters in $a\in M^l,$ then there is a uniformly $\emptyset $-definable family $\{X_t\}_{t\in T}$ of nonempty definable subsets of  $P$ such that $a\in T(M)$ and $X_{a}=X.$  We consider on $T$ the $\emptyset $-definable equivalence relation  given by $t\sim t'$ if and only if $X_t=X_{t'}$ and,  we denote by $[t]$ the $\sim $-equivalence class of $t\in T$ and $[T]=\{[t]:t\in T\}.$ In particular, below we consider $[a],$ a {\it code} for $X.$\\

\begin{prop}\label{prop def choice dcl}
Assume that ${\mathbb M}$ is sufficiently saturated. Then ${\mathbb M}$ has definable choice if and only if every nonempty definable subset $X\subseteq M$ defined with parameters $a\in M^l$ has an element in $\dcl([a]\,\mathbb{P})$.
\end{prop}

\pf
 ($\Rightarrow $) Let $X$ be a nonempty definable subset of $M$ defined with parameters $a\in M^l$. Then there is a uniformly $\emptyset $-definable family $\{X_t\}_{t\in T}$ of nonempty definable subsets of  $P$ such that $a\in T(M)$ and $X_{a}=X$. By Proposition \ref{prop def skolem/choice elementary substruture}, there is  a ${\mathbb P}$-definable choice function  $h:T\to P$ for the family $\{X_t\}_{t\in T}$ of the form $h=g(-,c)$ for some $\emptyset $-definable function $g:T\times D\to P$ and some $c\in P^q$. Since $g$ has the choice property, it induces a $\emptyset $-definable function $\bar{g}:[T] \times D\to M$. Then $\bar{g}([a],c)\in X\cap {\rm dcl}([a]\,{\mathbb P})$. 

($\Leftarrow $) Let $\{X_t\}_{t\in T}$ be a uniformly $\emptyset $-definable family of nonempty definable subsets of $M$. Suppose that $T\subseteq M^l$ and fix $v\in T$. Since $X_v$ is defined over $v$, there is an element  in   $X_v\cap {\rm dcl}([v]\,{\mathbb P})$.  So there is a $\emptyset$-definable function $g_v:Z_v\subseteq [T]\times M^{q_v}\to M$   and $c_v\in P^{q_v}$ such that $([v],c_v)\in Z_v$ and $g_v([v],c_v)\in X_v$. Let $B_v=\{[t] \in [T]: ([t],c_v)\in Z_v\}$. Then each $B_v$ is defined over ${\mathbb P}$.

Since $[T]=\bigcup _v B_v$, by saturation, there are finitely many $v_0,\ldots , v_m$ such that $[T]=B_{v_0}\cup \ldots \cup B_{v_m}$. For $i=0,\ldots, m$, let  $f_i: B_{v_i}\to M$ be given by $f_i([t])=g_{v_i}([t],c_{v_i})$. Then, by recursion on $m$, there is a ${\mathbb P}$-definable function $f:[T]\to M$ such that $f([t])\in X_t$ for all $t\in T$. Hence, the ${\mathbb P}$-definable function $h:T\to M$ given by $h(t)=f([t])$ is a definable choice function for the family $\{X_t\}_{t\in T}$. 
 \qed \\

\end{section}

\begin{section}{Preliminaries on the trichotomy theorem}\label{section prelim trichotomy}
\noindent
Recall that, in an o-minimal structure ${\mathbb M}=(M,<,\cdots)$, a point $a\in M$ is {\it non-trivial over $A\subseteq M$} if and only if there is an infinite open interval $I\subseteq M$ containing $a$ and a definable continuous function $F:I\times I\to M$ defined over $A$ such that $F$ is strictly monotone in each variable (separately). We say that $a\in M$ is {\it non-trivial} if it is non-trivial over some $A$. A point which is not non-trivial is called {\it trivial}.  Finally, let us say that an o-minimal structure ${\mathbb M}$ is {\it trivial} if and only if every $a\in M$ is trivial. 

In this section we recall the trichotomy theorem in an arbitrary o-minimal structure ${\mathbb M}$ and present several well known characterizations of trivial/non-trivial  that we require later.  

\medskip

First note the following (\cite[Lemma 2.2]{mrs}):

\begin{fact}\label{fact trivial J and K}
{\em
${\mathbb M}$ is non-trivial  if and only if there are infinite open intervals $J, K\subseteq M$ and a definable continuous function $G:J\times K\to M$ such that $G$ is strictly monotone in each variable.
}
\end{fact}

The main theorem of Peterzil and Starchenko in \cite{pst} is the following:

\begin{thm}[Trichotomy] Let ${\mathbb M} $ be a sufficiently saturated o-minimal structure structure. Given $a \in M$, then $a$ is non-trivial if and only if either
\begin{enumerate}
\item[(T2)] the structure that ${\mathbb M}$ induces on some convex neighbourhood of $a$ is an ordered vector space over an ordered division ring, or
\item[(T3)] the structure that ${\mathbb M}$ induces on some open interval around $a$ is an o-minimal expansion of a real closed field.
\end{enumerate}
\end{thm}

So, the Trichotomy theorem implies that if $a\in M$ is non-trivial, then $a$ is contained in a definable group-interval.  Recall that a {\it definable  group-interval}  $J=\langle (-b, b), 0, +,  <\rangle $ is an open interval $(-b, b)\subseteq M$, with $-b<b$ in $M\cup \{-\infty, +\infty \}$, together with a binary partial continuous definable operation $+:J^2\to J$ and an element $0\in J$, such that:
\begin{itemize}
\item[(i)]
 $x+y=y+x$ when defined;  
 $(x+y)+z=x+(y+z)$ when defined; 
 if  $x<y$ and $x+z$ and $y+z$ are defined then $x+z<y+z$;
\item[(ii)]
for every $x\in J$, if $x>0$, then the set $\{y\in J:\,\, x+y \,\, \textrm{is defined}\}$ is an interval of the form $(-b, r(x))$;
\item[(iii)]
for every $x\in J$, we have  $\lim _{z\to 0}(z+x)=x$ and if $x>0$ we also have   $\lim _{z\to r(x)^-}(x+z)=b$;
\item[(iv)]
for every $x\in J$ there exists $-x\in J$ such that $x+(-x)=0$. \\
\end{itemize}

Let us give other characterizations of trivial/non-trivial that will be useful later. First recall that an o-minimal structure ${\mathbb M}$ is a \emph{geometric structure}, i.e.:
\begin{itemize}
\item[-] 
$\mathrm{acl}(-)$ satisfies the Exchange Principle: if $a,b \in M$, $A \subseteq M$ and   $b \in \mathrm{acl}(aA)$, then either $b \in \mathrm{acl}(A)$ or $a \in \mathrm{acl}(bA)$ and,
\item[-] 
for any formula $\varphi(x,\bar{y})$, there exists $n \in \mathbb{N}$ such that for any $\bar{b}\in M^r$, either $\varphi(x,\bar{b})$ has fewer than $n$ solutions in ${\mathbb M} $ or it has infinitely many.
\end{itemize}
See \cite[Theorem 4.1]{pst}. So as in any geometric structure, given $A\subseteq M$, for $\bar{a}\in M^n$ we defined $\dim (\bar{a} / A)$ to be the minimal cardinality of $\bar{a}'\subseteq \bar{a}$ such that ${\rm acl }(\bar{a}' A)={\rm acl}(\bar{a} A)$, and for $U\subseteq M^n$ a definable subset over $A$, we define $\dim (U)=\max \{\dim (\bar{a} / A): \bar{a}\in U\}$. Furthermore, we say that $\bar{a}\in U$ is {\it generic in $U$ over $A$} if and only if $\dim (\bar{a} / A)=\dim(U)$. Note that in the o-minimal case, dimension of definable subsets over $A$ defined in this way coincides with the o-minimal dimension defined using the cell decomposition theorem (see for example \cite[Lemma 1.4]{Pi88} and \cite[Chapter 3, (2.7)]{vdd}).

Here is another characterization of trivial o-minimal structure related with the notion of trivial pregeometry (see \cite[Lemmas 1.8, 2.1 and 2.2]{mrs}). Recall that $\{a_1, \ldots, a_m\}\subseteq M$ is {\it independent} if there is no $1\leq i\leq m$ such that $a_i$ is algebraic over $\{a_1, \ldots, a_{i-1}, a_{i+1}, \ldots, a_m\}$.

\begin{fact}\label{fact trivial and independence}
${\mathbb M}$ is trivial  if and only if any finite subset of $M$ whose elements are pairwise independent is independent.
\end{fact}

Recall that  a \emph{curve} is a definable one-dimensional subset of $M^2$. A family of curves $\mathcal{F}$ is said to be \emph{definable} if there exist definable sets $U \subseteq M^k$ and $F\subseteq U \times M^2$ such that $\mathcal{F}=\{C_{\bar{u}}: \bar{u} \in U\}$ where for each $\bar{u}\in U$ we  let $C_{\bar{u}}:=\{\langle x, y \rangle: \langle \bar{u}, x, y \rangle \in F\}$. We say that $C_{\bar{u}}$ is \emph{generic} in $F$ if $\bar{u}$ is generic in $U$ over the parameters defining $F$. The family $F$ is said to be \emph{interpretable} if $U$ is replaced by $U/E$, where $E$ is a definable equivalence relation on $U$.
A definable (or interpretable) family of curves $\mathcal{F}=\{C_{\bar{u}}: \bar{u} \in U\}$ is \emph{normal} of dimension $n$ if $\dim(U)=n$ (or $\dim(U/E)=n$) 
and for $\bar{u} \neq \bar{v}$ from $U$ (or $U/E$), $C_{\bar{u}}$ and $C_{\bar{v}}$ intersect in at most finitely many points. \\

We have the following characterization given in \cite{pst}:

\begin{thm}\label{thm Z1 and curves}
Let ${\mathbb M} $ be a sufficiently saturated o-minimal structure. Then ${\mathbb M}$ is trivial if and only if the following holds:
\begin{enumerate}
\item[$(Z_1)$] for every interpretable infinite normal family of curves $\mathcal{F}$, if $\mathcal{C}$ is a generic curve in $\mathcal{F}$ and $\langle a,b\rangle$ is generic in $\mathcal{C}$, then either $\dim(\mathcal{C}\cap(\{a\} \times M))=1$ or $\dim(\mathcal{C}\cap(M \times \{b\}))=1$.
\end{enumerate}
\end{thm}

\medskip
We conclude with an observation that  is somewhat implicit in \cite{mrs} and follows from \cite[Chapter 2, (2.19) Exercises 1 and 2]{vdd}. First let us recall the material from the later source.  We say that an open cell $C\subseteq M^n$ is {\it regular} if for each $i\in \{1, \ldots, n\}$ and each two points $x, y\in C$ that differ only on the $i^{{\rm th}}$ coordinate and each point $z\in M^n$ that differs from $x$ and $y$ only in the  $i^{{\rm th}}$ coordinate, if $x_i<z_i<y_i$ then $z\in C$. 

We say that a definable function $f:C\subseteq M^n\to M$ on a regular open cell is {\it strictly increasing in the $i^{{\rm th}}$ coordinate} (resp. {\it strictly decreasing in the $i^{{\rm th}}$ coordinate} or {\it constant in the $i^{{\rm th}}$ coordinate}) if for each two points $x, y\in C$ that differ only on the $i^{{\rm th}}$ coordinate, if $x_i<y_i$ then $f(x)<f(y)$ (resp. $f(y)<f(x)$ and $(f(x)=f(y)$). 

We say that a definable function $f:C\subseteq M^n\to M$ on a regular open cell is {\it regular} if it is continuous, and for each $i\in \{1, \ldots, n\}$, $f$ is (depending on $i$) either strictly increasing or  strictly decreasing or  constant in the $i^{{\rm th}}$ coordinate.

Let us call a definable function $f:C\subseteq M^n\to M$ on a regular open cell {\it unary} if it is regular and is  constant in the $i^{{\rm th}}$ coordinate for all except for at most one $i\in \{1,\ldots, m\}$. \\

We have the following unary version of the Cell Decomposition Theorem for trivial o-minimal structures. We refer the reader to \cite[Chapter 3, (2.10)]{vdd} for the definition of a \emph{decomposition} of $M^n$.  

\begin{thm}[Unary CDT]\label{thm UCDT}
Suppose that ${\mathbb M}$ is a trivial o-minimal structure. Then
\begin{enumerate}
\item
Given any definable sets $A_1, \dots, A_k\subseteq  M^n$, there is a decomposition $\mathcal{C}$ of $M^n$ that partitions each $A_i$, all of whose open cells are regular.

\item
Given a definable function $f:A\subseteq M^n\to  M$, there is a decomposition $\mathcal{C}$ of $M^n$ that partitions $A$ all of whose open cells are regular, and such that for each open cell  $B\in\mathcal{C}$ with $B\subseteq  A$ the restriction $f_{|B}$  is unary.
\end{enumerate}
\end{thm}

\pf
\cite[Chapter 2, (2.19) Exercise 2]{vdd} shows that the regular cell decomposition theorem holds i.e., (1) above holds and (2) holds with the restrictions $f_{|B}$ being regular. Since  ${\mathbb M}$ is a trivial o-minimal structure, using Fact \ref{fact trivial J and K}, the restrictions $f_{|B}$ from (2) must be unary.
\qed \\

\end{section}

\begin{section}{Auxiliary lemmas}\label{section auxiliary lemmas}
\noindent
Here we let ${\mathbb M}$ be an arbitrary o-minimal structure and we prove a couple of auxiliary lemmas needed  in the next section. 

\begin{lem}\label{lem mu}
Let $Z\subseteq M$ be a finite union of open intervals defined over $s$. Suppose that we have a definable function $\mu : Z\to M$ defined over $s$ such that  $\mu (u)<u$ for all $u\in Z$. Then $Z$ is a finite union of points and open subintervals, all defined over $s$,  such that on each such open subinterval $X$, if  $u\in X$, then for all $v,w\in (\mu (u), u)$ with  $v<w$ we have $v\in (\mu (w),w)$. 
\end{lem}

\pf
By the monotonicity theorem (\cite[Chapter 3, (1.2)]{vdd}), $Z$ is a finite union of points and open subintervals such that on each open subinterval $X$, the restriction $\mu :X\to Y_s$ is either constant or strictly monotone and continuous. Let $u\in X$ and suppose that $v,w\in (\mu (u), u)$ with  $v<w$.

Observe that  $\mu (u)<v<w$. It follows that: 
If $\mu $ is constant, then $\mu (u)=\mu (w)=\mu (v)$ and so   $v\in (\mu (w),w)$. 
If $\mu $ is strictly increasing, then  $\mu (v)<\mu (w)<\mu (u)$ and therefore $v\in (\mu (w),w)$.  
If $\mu $ is strictly decreasing, then  $\mu (w)<\mu (v)$ and hence  $v\in (\mu (w),w)$. 
\qed \\

Similarly:

\begin{lem}\label{lem mu+}
Let $Z\subseteq M$ be a finite union of open intervals defined over $s$. Suppose that we have a definable function $\mu : Z\to M$ defined over $s$ such that  $u<\mu (u)$ for all $u\in Z$. Then $Z$ is a finite union of points and open subintervals, all defined over $s$,  such that on each such open subinterval $X$, if  $u\in X$, then for all $v,w\in (u, \mu (u))$ with  $v<w$ we have $w\in (v, \mu (v))$. 
\end{lem}



\medskip

The following lemma is obtained combining Lemma \ref{lem mu} and the Unary CDT (Theorem \ref{thm UCDT}):

\begin{lem}\label{lem trivial and dim 2 family}
Let $Z\subseteq M$ be a finite union of  open intervals and $b :Z\to M$ a definable function defined over $s$ such that  $b(u)<u$ for all $u\in Z$. Let $\Omega : (b, \id _Z)_Z\to M$ be a definable function over $s$ such that for each $u\in Z$ we have:
\begin{itemize}
\item[-]
$\Omega (u, -)$ is strictly monotone and continuous.
\end{itemize}
Suppose that $Z$ has no non-trivial points. Then $Z$ is a finite union of points and open subintervals, all defined over $s$,  such that on each such open subinterval $X$:
\begin{itemize}
\item[-]
there is a definable function $d_X :X\to M$ defined over $s$  such that  $b_{|X}(u)\leq d_X(u)<u$ for all $u\in X$;
\item[-]
for every  $u\in X$, if $v, w\in (d_X(u), u)\cap X$ are such that  $v<w$,  then we have $v\in (d_X(w), w)$ and $\Omega (u,v)=\Omega (w,v)$.
\end{itemize}
\end{lem}

\pf
By Theorem \ref{thm UCDT}, take $\mathcal{D}$ a regular cell decomposition of $M^2$ that partitions $C:=(b, \id _Z)_Z$ such that for each open cell $D\in \mathcal{D}$ contained in $C$ the restriction $\Omega _{|D}$ is unary.  

Let $Z_1, \ldots, Z_k$ be the open cells of $\pi (\mathcal{D})$ contained in $Z$.  For each $i$ let $D_i$ be the open cell of $\mathcal{D}$ of the form $(d_i, \id _{Z_i})_{Z_i}$ where $d_i:Z_i\to M$ is the maximal element of the set of continuous definable function with domain $Z_i$, which are below $\id _{Z_i}$ and whose graph is a cell of $\mathcal{D}$. Of course $b_{|Z_i}\leq d_i$.

Fix $i$. Since for all $u\in Z_i$ we have $d_i(u)<u$, by Lemma \ref{lem mu},  $Z_i$ is a finite union of points and open subintervals $Z_{i,1},\ldots, Z_{i,l_i}$, all defined over $s$,  such that for each $j\in \{1, \ldots, l_i\}$, if  $u\in Z_{i,j}$, then for all $v,w\in (d_i(u), u)$ with  $v<w$ we have $v\in (d_i(w),w)$.   

Let $u\in Z_{i,j}$ and let $v, w\in (d_i(u), u)\cap Z_{i,j}$ with $v<w$. Then we have $v\in (d_i(w), w)$ and therefore $(u,v), (u,w), (w,v)\in D_i$. Since $D_i\in \mathcal{D}$ and $\Omega _{|D}$ is unary we have $\Omega (u,v)=\Omega (w,v)$ as required.
\qed \\

Similarly, using Lemma \ref{lem mu+} instead of Lemma \ref{lem mu} we have:

\begin{lem}\label{lem trivial and dim 2 family+}
Let $Z\subseteq M$ be a finite union of  open intervals and $b :Z\to M$ a definable function defined over $s$ such that  $u<b(u)$ for all $u\in Z$. Let $\Omega : (\id _Z, b)_Z\to M$ be a definable function over $s$ such that for each $u\in Z$ we have:
\begin{itemize}
\item[-]
$\Omega (u, -)$ is strictly monotone and continuous.
\end{itemize}
Suppose that $Z$ has no non-trivial points. Then $Z$ is a finite union of points and open subintervals, all defined over $s$,  such that on each such open subinterval $X$:
\begin{itemize}
\item[-]
there is a definable function $d_X :X\to M$ defined over $s$  such that  $u<d_{X}(u)\leq b_{|X}(u)$ for all $u\in X$;
\item[-]
for every  $u\in X$, if $v, w\in (u, d_X(u))\cap X$ are such that  $v<w$,  then we have $w\in (v, d_X(v))$ and $\Omega (u,w)=\Omega (v,w)$.
\end{itemize}
\end{lem}

\medskip
We end with the following lemma which will be used later to ``remove'' parameters from certain uniformly definable families:

\begin{lem}\label{lem removing parameters}
Assume that $\mathbb{M}$ has definable Skolem functions. Consider the following data $(Z, T, \{Z_t\}_{t\in T}, \{\mathcal{U}_t\}_{t\in T})$ where: $Z\subseteq M$ is a $\emptyset $-definable open subinterval, $T\subseteq Z$ is a $\emptyset $-definable subset such that $Z\setminus T$ is finite,  $\{Z_t\}_{t\in T}$  is a uniformly $\emptyset $-definable family of open subintervals of $Z$ such that $Z\setminus \bigcup Z_t$ is finite and $\{\mathcal{U}_t\}_{t\in T}$  is a uniformly $\emptyset $-definable family of finitely many two by two disjoint open subintervals of each $Z_t$ such that $Z_t\setminus \bigcup \mathcal{U}_t$ is finite. Then after naming finitely many elements of $\mathbb{P}$, and after removing  finitely many $\emptyset $-definable points from $Z$ and replacing $Z$ by one of the remaining open subintervals, we may assume that there are $\emptyset$-definable functions $q^-, q^+:Z\to M$  such that $q^-<q^+$ and $t:Z\to T$ such that, for each $z\in Z$,  $(q^-(z), q^+(z))$ is one of the open subintervals of $Z_{t(z)}$ in $\mathcal{U}_{t(z)}$ and $z\in (q^-(z), q^+(z))$.
\end{lem}

\pf
First we show that 
\begin{clm}
The set $Z':=\{z\in Z: z\notin \bigcup \mathcal{U}_t\,\,\textrm{for all}\,\,t\in T\}$ is finite and defined over $\emptyset $. 
\end{clm}

\pf
Suppose that $Z'$ is infinite. By the monotonicity theorem and the  finiteness theorem (see \cite[Chapter 3, (1.8)]{vdd}) there is a partition of $T$ into  finitely many points and open intervals, all defined over $\emptyset $, such that on each such open subinterval $J$, there are $n_J$ continuous definable functions $f_{J, i}:J\to M\cup \{-\infty, +\infty \}$ such that $-\infty \leq f_{J,1}<\ldots <f_{J, n_J}\leq +\infty $  and for each $t\in J$ the end points of the open subintervals in $\mathcal{U}_t$ are exactly  $f_{J,1}(t)<\ldots <f_{J, n_J}(t)$.

Fix such $J$ and let $Z'_J:=\{z\in Z: z\notin \bigcup \mathcal{U}_t\,\,\textrm{for all}\,\, t\in J\}$. Since $Z'\subseteq Z'_J$, we have that $Z'_J$ is infinite. On the other hand,  for each $z\in Z'_J$ and each $t\in J$, there is a unique index $i(z,t)$ such that $z=f_{J, i(z,t)}(t)$.  So we have a definable function $i:Z'_J\times J\to \{1, \ldots, n_J\}$ and, by the cell decomposition theorem, there is a finite partition $\mathcal{C}$ of $Z'_J\times J$ into finitely many cells such that on each such cell $C$ the function is constant, say with value $i_{C}$. 

Since $J$ is an open subinterval and $Z'_J$ is infinite, there is an open cell $C\in \mathcal{C}$ of the form $(f,g)_B$ with $B\subseteq Z'_J$ an open subinterval and $f,g: B\to M\cup \{-\infty, +\infty \}$ continuous definable functions such that $-\infty \leq f<g\leq +\infty $. Fix such a cell $C$ and note that for all $z\in B$  the definable function $f_{J,i_C}: (f(z), g(z))\to M\cup \{-\infty , +\infty \}$ is constant with value $z$. 

Now fix $z'\in B$ and let $I$ be a closed infinite interval contained in $(f(z'), g(z'))$. Then by continuity there is an open subinterval $O \subseteq B$ such that  for all $z\in O$ we have $I\subseteq (f(z), g(z))$. This implies that for each $z\in O$, the definable function $f_{J,i_C}: I\to M\cup \{-\infty , +\infty \}$ is constant with value $z$, which is absurd since $O$ is infinite.
\qed \\

It follows that for each $z\in Z\setminus Z'$, the set $T_z=\{t\in T: z\in \bigcup \mathcal{U}_t\}$ is a nonempty definable set. So $\{T_z\}_{z\in Z\setminus Z'}$ is a uniformly definable family of nonempty definable sets, defined over $\emptyset $. By definable Skolem functions, there is a definable function $t: Z\setminus Z'\to T$ such that for each $z\in Z\setminus Z'$, we have $z\in \bigcup \mathcal{U}_{t(z)}$. Varying the parameter over which $t:Z\setminus Z'\to T$ is defined, we get a uniformly $\emptyset $-definable family of such functions. By Proposition \ref{prop def skolem dcl} we can choose a parameter in ${\rm dcl}(p)$ for some $p\in P^k$. After naming the elements in $p$ we see that we assume that $t: Z\setminus Z'\to T$  is defined over $\emptyset $.

Now $\{Z_{t(z)}\setminus \mathcal{U}_{t(z)}\}_{z\in Z\setminus Z'}$ is a uniformly $\emptyset $-definable family of finite sets. By Finiteness lemma (\cite[Chapter 3 (1.7)]{vdd}), the cardinality of such finite sets are uniformly bounded. Hence there is a finite definable subset $Z''\subseteq Z$ such that $Z\setminus (Z'\cup Z'')$ is a finite union of open subintervals, all defined over $\emptyset$, such that on each such open subinterval the cardinality of all the finite sets is fixed. So replacing $Z$ by one of the remaining open subintervals of $Z\setminus (Z'\cup Z'')$,   we can assume that there is $n_Z$ such that $\# (Z_{t(z)}\setminus \mathcal{U}_{t(z)})=n_Z$ for every $z\in Z$.   

Consider the $\emptyset$-definable functions $q_i:Z\to M$ ($i=1, \ldots, n_Z$) such that $q_{1}<\ldots <q_{n_Z}$ and, for each $z\in Z$,  $(q_{i}(z), q_{i+1}(z))$ are the finitely many open subintervals of $Z_{t(z)}$ in $\mathcal{U}_{t(z)}$. Since the open subintervals of $\mathcal{U}_{t(z)}$ are two by two disjoint, the subsets $\{z\in Z: z\in (q_{i}(z), q_{i+1}(z))\}$ with $i=1, \ldots , n_Z$, determine a partition of $Z$ into finitely many $\emptyset $-definable subsets. These subsets are each a finite union of open intervals and points all defined over $\emptyset $, so removing all those finitely many points and replacing $Z$ by one of the remaining open subintervals, we may assume that there are $\emptyset$-definable functions $q^-, q^+:Z\to M$  such that $q^-<q^+$ and, for each $z\in Z$,  $(q^-(z), q^+(z))$ is one of the open subintervals of $Z_{t(z)}$ in $\mathcal{U}_{t(z)}$ and $z\in (q^-(z), q^+(z))$.
\qed \\

\end{section}

\begin{section}{Technical lemmas}\label{section technical lemmas}

{\it Now we assume that ${\mathbb M}$ is a sufficiently saturated o-minimal structure with  definable Skolem functions}.  \\

First, using existence of definable Skolem functions,  we try to find finitely many open subintervals $I$ of $M$ obtained by removing from $M$ finitely many points, such that on each one of them we have a definable function $H:I\times I\to M$ which is strictly monotone in each variable. This goal will only be partially achieved in Lemmas \ref{lem def skolem and family of funct} and \ref{lem H_s}. We will get instead on each such subinterval a uniformly definable family $\{H_s:Y_s\times I_s\to M\}_{s\in S}$ of  functions which are strictly monotone in  the second variable (when the first variable is fixed). \\

For each $u\in M$  let 
$$L_u=\{(x,y)\in M^2: x<y<u\}$$
and let
$$R_u=\{(x,y)\in M^2: u<y<x\}.$$

\begin{lem}\label{lem def skolem and family of funct}
After naming finitely many elements of $\mathbb{P}$, there is a finite collection $\mathcal{I}$ of $\emptyset$-definable open intervals of $M$ such that:
\begin{itemize}
\item[(1)] 
$M\setminus \bigcup \mathcal{I}$ is finite.
\item[(2)]
For each $I\in {\mathcal I}$ there are $\emptyset$-definable functions $l, r:I\to M$, with $l<r$, which are either constant or strictly monotone and continuous and, there is a uniformly $\emptyset $-definable family
$$\{h_{I,u}:(l(u),r(u))\to M\}_{u\in I}$$ 
of definable, continuous and  strictly increasing  functions such that for each $u\in I:$
\begin{itemize}
\item[-]
$(t,h_{I,u}(t))\in L_u$ for all $l(u)<t<u$ and
\item[-] 
$(t, h_{I,u}(t))\in R_u$ for all $u<t<r(u)$.
\end{itemize}
\end{itemize}
\end{lem}

\pf
Fix $u\in M$. By definable Skolem functions there is a definable function $f_u:(-\infty ,u)\to (-\infty ,u)$ such that for each $t\in (-\infty ,u)$ we have $(t,f_u(t))\in L_u$  and a definable function $g_u:(u,+\infty )\to (u,+\infty )$ such that for each $t\in (u,+\infty )$ we have $(t,g_u(t))\in R_u$. Varying the parameter over which $f_u$ (resp. $g_u$) is defined, we get a uniformly definable family of such functions defined over $u$. By Proposition \ref{prop def skolem dcl} we can choose a parameter in ${\rm dcl}(u\,p)$ for some $p\in P^k$. After naming the elements in $p$ we see that we assume that $f_u$  (resp. $g_u$) is defined over $u$.  

By the monotonicity theorem (\cite[Chapter 3, (1.2)]{vdd}), there is $l_u\in (-\infty ,u)$ such that $f_{u|(l_u,u)}$ is either constant, or strictly monotone and continuous and, there is  $r_u\in (u,+\infty )$ such that $g_{u|(u,r_u)}$ is either constant or strictly monotone.  Note that $f_{u}$ on $(l_u,u)$ cannot be constant, for if $f_u(t)=c$ for all $t\in (l_u,u),$ then $l_u<c<u$ and so for $s\in (c,u)$ we would have $c<s<f_u(s)=c.$ Note also that $f_u$ on $(l_u,u)$ cannot be strictly decreasing, otherwise, taking $d\in (l_u,u),$ then for all $t\in (d,u)$ we would have $f_u(t)<f_u(d)<u,$ consequently for $s\in (f_u(d),u)\subseteq (d,u)$ we would have $f_u(d)<s<f_u(s)<f_u(d).$ So $f_u$ must be strictly increasing and continuous on $(l_u,u)$. Similarly, $g_u$ must be strictly increasing and continuous on $(u, r_u)$.


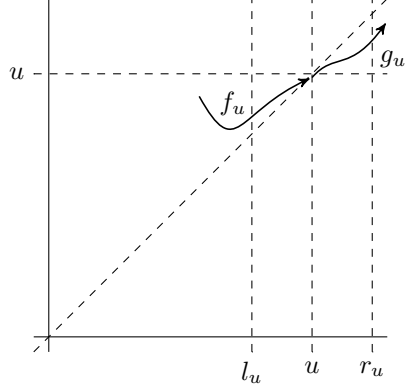
\begin{figure}[ht]
 \caption{The functions $f_u$ and $g_u$}
\begin{center}
\begin{tikzpicture}\label{fig fu and gu}
  \draw[->,>=stealth',shorten >=1pt,auto,node distance=2.8cm,
                    semithick] (3.5,3.45) .. controls (3.8,3.8) and (4.0,3.5) ..  node[below, right] {$\,\,\,\,\,\,g_u$}  (4.5,4.2);
     \draw[->,>=stealth',shorten >=1pt,auto,node distance=2.8cm,
                    semithick] (2.0,3.2) .. controls (2.5,2.3) and (2.5,3.0) ..  node[above] {$f_u\,\,\,\,$}  (3.5,3.45);
     \draw[dashed] (-0.2,-0.2) -- (4.5,4.5);
     \draw[dashed] (3.5,-0.2) -- node[below]{$u$}(3.5,-0.2);
      \draw[dashed] (3.5,-0.2) -- (3.5,4.5);
     \draw[dashed] (-0.2,3.5) -- node[left]{$u$}(-0.2,3.5);
     \draw[dashed] (-0.2,3.5) -- (4.5,3.5);
     \draw[dashed] (2.7,-0.2) -- node[below]{$l_u$}(2.7,-0.2);   
      \draw[dashed] (2.7,-0.2) -- (2.7,4.5); 
     \draw[dashed] (4.3,-0.2) -- node[below]{$r_u$}(4.3,-0.2);  
      \draw[dashed] (4.3,-0.2) -- (4.3,4.5);
     \draw[] (-0.2,0.0) -- node[below]{$$}(4.5,0.0); 
     \draw[] (0.0,-0.2) -- node[below]{$$}(0.0,4.5); 
\end{tikzpicture}
\end{center}
\end{figure}

Note also that $\lim _{t\to u^-}f_u(t)=u=\lim _{t\to u^+}g_u(t)$. Let $h_u:(l_u, r_u)\to M$ given by
\begin{eqnarray*}
h_u(t)=
\begin{cases}
f_u(t) \qquad \textrm{if}\quad t\in (l_u, u)\\
\\
u \qquad  \quad \,\,\,\,\,\textrm{if}\quad t =u\\
\\
g_u(t)\qquad  \,\,\textrm{if}\quad t\in (u, r_u)
\end{cases}
\end{eqnarray*}
Then, for each $u\in M,$ there is a definable function $h_u:(l_u, r_u)\to M$  such that $h_u$ is continuous and strictly increasing, $(t,h_u(t))\in L_u$ for all $t\in (l_u,u)$, $(t, h_u(t))\in R_u$ for all $t\in (u, r_u)$ and $\lim _{t\to u^-}h_u(t)=u=\lim _{t\to u^+}h_u(t)$. Furthermore,  $h_u$ is defined over $u$.

Consider the collection of all first-order formulas $\varphi (x,y,v)$ over $\emptyset $ such that for some $a\in M$, $\varphi (x,y,a)$ defines the graph of a function on some open interval containing $a$ with the same properties as $h_a$ above.  Since all of the above properties of $h_u$ are first-order on $u$, for each such $\varphi $, the  subset $I_{\varphi }$ of $M$ of all $a\in M$ such that $\varphi (x,y,a)$ defines the graph of a function on some open interval containing $a$ with the same properties as $h_a$ above is  $\emptyset $-definable. Let  ${\mathcal I}$ be the collection of all those $\emptyset$-definable subsets of $M$. Then $\bigcup {\mathcal I}=M$ and for each $I\in {\mathcal I}$ there is a uniformly $\emptyset $-definable family
$$\{h_{I,u}:(l_u, r_u)\to M\}_{u\in I}$$ 
of definable, continuous and  strictly increasing  functions such that for each $u\in I:$
\begin{itemize}
\item[-]
$(t,h_{I,u}(t))\in L_u$ for all $t\in (l_u, u)$,
\item[-]
$(t, h_{I,u}(t))\in R_u$ for all $t\in (u, r_u)$ and,
\item[-]
$\lim _{t\to u^-}h_{I,u}(t)=u=\lim _{t\to u^+}h_{I,u}(t)$.
\end{itemize}
By saturation we may assume that ${\mathcal I}$ is finite. By o-minimality, we may also assume,  after ignoring finitely many elements of $M$, that each $I\in {\mathcal I}$ is a $\emptyset $-definable open interval of $M$.

Fix $I\in {\mathcal I}$ and consider the $\emptyset$-definable functions $l, r:I\to M$ given by $l(v)=l_v$ and $r(v)=r_v$.  By the monotonicity theorem (\cite[Chapter 3, (1.2)]{vdd}), after removing  finitely many points from each $I$ we can  assume that on each $I\in {\mathcal I}$ the $\emptyset$-definable functions $l, r:I\to M$ are  constant or strictly monotone and continuous. 
\qed \\

Given $Z\subseteq M$ let us set $\Delta _Z=\{(x,x): x\in Z\}$.

\begin{lem}\label{lem H_s}
Let $\mathcal{I}$ and for each $I\in \mathcal{I}$ let $\{h_{I,w}: (l(w), r(w))\to M\}_{w\in I}$ be as in Lemma \ref{lem def skolem and family of funct}. Then after naming more finitely elements of $\mathbb{P}$, for each $I\in \mathcal{I}$ there is a finite collection $\mathcal{Y}^I$ of $\emptyset$-definable open subintervals of $I$ such that:
\begin{itemize}
\item[-]
$\bigcup \mathcal{Y}^I=I$.
\item[-]
For each $Y\in \mathcal{Y}^I$ there is a $\emptyset $-definable $e_Y\in Y$ and a uniformly $\emptyset$-definable family of definable functions 
$$\{H_s:Y_s\times I_s\to M\}_{s\in S}$$
 given by $H_s(w,t)=h_{I,w}(t)$ where, $S=Y\setminus \{e_Y\}$, for each $s\in S$,  
\begin{eqnarray*}
Y_s=
\begin{cases}
Y\cap (s,e_Y)\qquad \textrm{if}\,\,\, s<e_Y\\
\\
Y\cap (e_Y,s)\qquad \textrm{otherwise}
\end{cases}
\end{eqnarray*}
and $I_s=(l_s,r_s)\subseteq I$ with $l_s=\max \{l(s), l(e_Y)\}$ and $r_s=\min \{r(s), r(e_Y)\}$. 
Moreover, $\Delta _{Y_s}\subseteq Y_s\times I_s$.
\end{itemize}
\end{lem}

\pf
For each $u\in I$, since $l(u)<u<r(u)$, by continuity, there is an open interval $Y_u\subseteq I$ containing $u$ such that for all $v\in Y_u$ we have $l(v)<u<r(v)$.  Since the property of $Y_u$ is first-order on $u$, there is a collection $\mathcal{Y}^I$ of $\emptyset$-definable subsets of $I$  such that $\bigcup \mathcal{Y}^I=I$ and for each $Y\in \mathcal{Y}^I$, $Y$ is a $\emptyset$-definable open subinterval of $I$ and for all $x, x'\in Y$ we have $l(x')<x<r(x')$. By saturation, $\mathcal{Y}^I$ is finite.

Fix $Y\in \mathcal{Y}^I$. By Proposition \ref{prop def skolem dcl}, we can choose an element of $Y$ in ${\rm dcl}(p)$ for some $p\in P^k$. After naming the elements in $p$ we see that we can  choose a $\emptyset$-definable $e_Y\in Y$. 

For $s\in Y\setminus \{e_Y\}$, let $l_s$, $r_s$, $I_s$ and  $Y_s$ be as in the statement of the lemma. Let $w\in Y_s\subseteq Y\subseteq I$. Since the $\emptyset$-definable functions $l, r:I\to M$ are  constant or strictly monotone and continuous and $w$ is between $s$ and $e_Y$, we have $(l_s, r_s)\subseteq (l(w), r(w))$. Moreover, since $s, e_Y\in Y$ we also have $l_s<w <r_s$. It follows that  $\Delta _{Y_s}\subseteq Y_s\times I_s$ and we can define uniformly the definable functions 
$$H_s:Y_s\times I_s\to M$$
by $H_s(w,t)=h_{I,w}(t)$.
\qed \\

Each $H_s$ of Lemma \ref{lem H_s} is  continuous and strictly increasing in the second coordinate when we fix the first coordinate. We need to  investigate what happens in the first coordinate when we fix the second coordinate. Following this line of thought does not seem to produce any conclusive results, so we try to see what happens to the restrictions of the uniformly definable family of curves
$$C_{w}=\{(t,f_w(t)): t\in I_{s,w}^-\},$$
with  $w\in Y_{s}$ and $I_{s,w}^-=\{t\in I_s:t<w\}$, induced by $\{H_s\}_{s\in S}$.

In Lemmas \ref{lem type (N1) or type (L)} and \ref{lem type (L)} below, we will find finitely many definable open subintervals of each $Y_s$, obtained by removing finitely many points, uniformly in $s$,  such that on each one of them the family $\{C_w\}_{w\in Y_s}$ induces one-dimensional definable normal families of curves or if not then we have density of non-trivial points on the subinterval.  \\

Below for each  $u\in Y_s$  we let $Y_{s,u}^-=\{v\in Y_s: v<u\}$.

\begin{lem}\label{lem type (N1) or type (L)}
For each $Y_s$ there is a finite collection $\mathcal{X}_s$,  uniformly defined over $s$, of two by two disjoint  open intervals uniformly defined over $s$ such that $Y_s\setminus \bigcup \mathcal{X}_s$ is finite and for each $X\in \mathcal{X}_s$ either  the following holds\\
\begin{itemize}
\item[(N1)] 
for every $u\in X$ there is an open subinterval $X_{s,u}^-\subseteq Y_{s,u}^-$  with right endpoint $u$  such that for every  $v\in X_{s,u}^-$ there is an open subinterval $J_{s,u,v}^-\subseteq I_{s,v}^-$ with right endpoint $v$ such that $\{C_{v|J_{s,u,v}^-}\}_{v\in X_{s,u}^-}$ is a one-dimensional definable  normal family of  curves,\\
\end{itemize}
or if not, then the following holds
\begin{itemize}
\item[(L)] 
for every $u\in X$ there is an open subinterval $X_{s,u}^-\subseteq Y_{s,u}^-$  with right endpoint $u$  such that  for every $v\in X_{s,u}^-$ we have $\bar{C}_{v|J_{s,u,v}^-}\cap C_{u|J_{s,u}^-}$ is finite where $J_{s,u,v}^-$ (resp. $J_{s,u}^-$) is an open subinterval of $I_{s,v}^-$ (resp. $I_{s,u}^-$) with right endpoint $v$ (resp. $u$),
\end{itemize}
or else
\begin{itemize}
\item[(D1)] 
every open subinterval of $X$ defined over $s$ has a non-trivial point.
\end{itemize}
\end{lem}

\proof
We are going to construct  definable subsets 
$$\mathcal{T}_s\subseteq \mathcal{B}_s\subseteq Y_s,$$
all defined over $s$, so that
$$Y_s=\mathcal{T}_s\sqcup  (\mathcal{B}_s\setminus \mathcal{T}_s)\sqcup (Y_s\setminus \mathcal{B}_s).$$
Each one of these pieces of $Y_s$ is a finite union of points and open intervals. We let $\mathcal{X}_s$ be the collection of all such open subintervals. Let $X\in \mathcal{X}_s$. We show that: 
\begin{itemize}
\item[-]
if $X$ is an open subinterval  of $Y_s\setminus \mathcal{B}_s$ then  (N1)  of the lemma holds on $X$;
\item[-]
if $X$ is an open subinterval of  $\mathcal{B}_s\setminus \mathcal{T}_s$  then   (L)  of the lemma holds on $X$;
\item[-]
if $X$ is an open subinterval of $\mathcal{T}_s$  then  (D1)  of the lemma holds on $X$.
\end{itemize}

Fix $u\in Y_s$. Consider the definable subset
$$\mathcal{B}_{s,u}=\{v\in Y_{s,u}^-: \,\textrm{for some}\,\,t<v\,\,\textrm{we have}\,\,f_u(t)\leq f_v(t)\}$$ 
and the definable function $\tau _{s,u}:\mathcal{B}_{s,u}\to I_{s,u}^-$ given by 
$$\tau _{s,u}(v)=\sup \{t<v:f_u(t)\leq f_v(t)\}.$$ 

Note that the family $\{\tau _{s,u}:\mathcal{B}_{s,u}\to I_{s,u}^-\}$ is uniformly defined over $s$ and $u$.

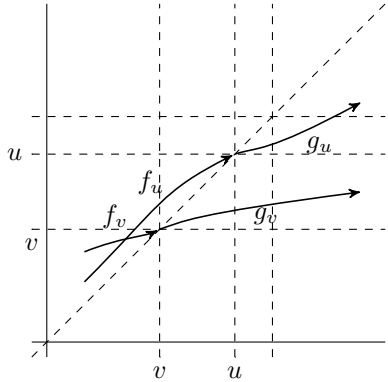
\begin{figure}[ht]
 \caption{The functions $f_u$, $g_u$, $f_v$ and $g_v$}
\begin{center}
\begin{tikzpicture}\label{fig fu and gu}
  \draw[->,>=stealth',shorten >=1pt,auto,node distance=2.8cm,
                    semithick] (2.5,2.5) .. controls (2.8,2.6) and (2.9,2.5) ..  node[below, right] {$\,\,\,\,\,\,g_u$}  (4.2,3.2);
     \draw[->,>=stealth',shorten >=1pt,auto,node distance=2.8cm,
                    semithick] (0.5,0.8) .. controls (1.5,1.8) and (1.5,2.0) ..  node[above] {$f_u\,\,\,\,$}  (2.5,2.5);
    \draw[->,>=stealth',shorten >=1pt,auto,node distance=2.8cm,
                    semithick] (1.5,1.5) .. controls (1.8,1.6) and (1.9,1.7) ..  node[below, right] {$\,\,\,\,\,\,\,\,\,g_v$}  (4.2,2.0);
     \draw[->,>=stealth',shorten >=1pt,auto,node distance=2.8cm,
                    semithick] (0.5,1.2) .. controls (1.0,1.4) and (1.3,1.4) ..  node[above] {$f_v\,\,\,\,\,\,\,$}  (1.5,1.5);                
     \draw[dashed] (-0.2,-0.2) -- (4.5,4.5);
     \draw[dashed] (2.5,-0.2) -- node[below]{$u$}(2.5,-0.2);
      \draw[dashed] (2.5,-0.2) -- (2.5,4.5);
     \draw[dashed] (-0.2,2.5) -- node[left]{$u$}(-0.2,2.5);
     \draw[dashed] (-0.2,2.5) -- (4.5,2.5);
     \draw[dashed] (-0.2,1.5) -- node[below]{$v$}(-0.2,1.5);  
     \draw[dashed] (-0.2,1.5) -- (4.5,1.5);  
     \draw[dashed] (1.5,-0.2) -- node[below]{$v$}(1.5,-0.2);  
     \draw[dashed] (1.5,-0.2) -- (1.5,4.5); 
      \draw[dashed] (3.0,-0.2) -- (3.0,4.5);
     \draw[dashed] (-0.2,3.0) -- (4.5,3.0);  
     \draw[] (-0.2,0.0) -- node[below]{$$}(4.5,0.0); 
     \draw[] (0.0,-0.2) -- node[below]{$$}(0.0,4.5); 
\end{tikzpicture}
\end{center}
\end{figure}

\begin{clm}\label{clm t_v<v}
For all $v\in \mathcal{B}_{s,u}$, we have $\tau _{s,u}(v)<v$ and $f_u(\tau _{s,u}(v))=f_v(\tau _{s,u}(v))$. 
\end{clm}

\pf
Note that for all $v<u$, since $f_u(v)\in L_u$, we have $v<f_u(v)$. So by continuity there is an open definable neighbourhood $V$ of $v$ such that for all $t\in V$ we have $v<f_u(t)<u$. In particular, there is $t'<v$ an element of $V$ such that $\{t<v: f_u(t)\leq f_v(t)\}$ is below $t'$, since for $t<v$ we have $f_v(t)<v$. It follows  that $\tau _{s,u}(v)<v$.

Suppose that $f_v(\tau _{s,u} (v))<f_u(\tau _{s,u}(v))$. Then by continuity there is an open definable neighbourhood $U$ of $\tau _{s,u}(v)$ such that for all $t\in U$ with $t<\tau _{s,u}(v)$ we have $f_v(t)<f_u(t)$, contradicting the definition of $\tau _{s,u}(v)$. Similarly, if  $f_u(\tau _{s,u}(v))< f_v(\tau _{s,u}(v))$, then by continuity there is an open definable neighbourhood $U$ of $\tau _{s,u}(v)$ such that for all $t\in U$ with $\tau _{s,u}(v)<t<v$ we have $f_u(t)<f_v(t)$, contradicting again the definition of $\tau _{s,u}(v)$. Therefore, $f_u(\tau _{s,u}(v))=f_v(\tau _{s,u}(v))$.

\qed \\

For each $u\in Y_s$, $\mathcal{B}_{s,u}$ is a finite union of points and open subintervals of $Y_{s,u}^-$. The collection of these finitely many points and open subintervals is uniformly defined over $s$ and $u$. Therefore, either $\mathcal{B}_{s,u}$ or its complement  in $Y_{s,u}^-$, but not both, contains an open subinterval with right endpoint $u$. Let
$$\mathcal{B}_s=\{u\in Y_s: \mathcal{B}_{s,u}\,\,\textrm{contains an open subinterval with right endpoint}\,\, u\}.$$ 

Note that $\mathcal{B}_s$ is uniformly defined over $s$.

\medskip
\noindent
{\bf Inside $Y_s\setminus \mathcal{B}_s$.} We have that $Y_s\setminus \mathcal{B}_s$ is a finite union of points and open subintervals of $Y_s$, all uniformly defined over $s$. Removing these finitely many points we may assume that $Y_s\setminus \mathcal{B}_s$ is a finite union of open subintervals of $Y_s$  uniformly  defined over $s$.

Let $u\in Y_s\setminus \mathcal{B}_s$. Then there is a maximal open subinterval $X_{s,u}^-$ with right endpoint $u$ such that $\mathcal{B}_{s,u}\cap X_{s,u}^-=\emptyset $. 

Then for all $v\in X_{s,u}^-$ we have  $f_v(t)<f_u(t)$ for all $t\in I_{s,v}^-$. It follows that,  for all $v\in X_{s,u}^-$ we have $C_{v}\cap C_{u}=\emptyset $.  

Since $\{X_{s,u}^-\}_{u\in Y_s\setminus \mathcal{B}_s}$ is uniformly definable, we have a definable map $\mu :Y_s\setminus \mathcal{B}_s\to Y_s$ given by $\mu (u)=\inf X_{s,u}^-$. So we have $X_{s,u}^-=(\mu (u),u)$ and  $\mu (u)<u$ for all $u\in Y_s\setminus \mathcal{B}_s$. 

Now by Lemma \ref{lem mu}, the set  $Y_s\setminus \mathcal{B}_s$ is a finite union of open subintervals $X$ all defined over $s$ such that, if  $u\in X$ and  $v,w\in X_{s,u}^-$ are such that $v<w$, then $v\in X_{s,w}^-$. But by definition of $X_{s,w}^-$, we have  $C_v\cap C_w=\emptyset $. Hence $\{C_v\}_{v\in X_{s,u}^-}$ is a one-dimensional definable normal family of curves, and  (N1) holds on $Y_s\setminus \mathcal{B}_s$.

\medskip
We have that $\mathcal{B}_s \subseteq Y_{s}$ is a finite union of points and open subintervals of $Y_s$, all uniformly defined over $s$. Removing these finitely many points we may assume that $\mathcal{B}_s$ is a finite union of open subintervals of $Y_s$  uniformly defined over $s$.

Let $u\in \mathcal{B}_s$. Then $\mathcal{B}_{s,u}$ contains a maximal open subinterval with right endpoint $u$. By abuse of notation we assume that $\mathcal{B}_{s,u}$ is that subinterval. These open subinterval are uniformly defined over $s$. Let
$$\mathcal{T}_s=\{u\in \mathcal{B}_s: \sup \tau _{s,u}=u\}.$$

Note that $\mathcal{T}_s$ is uniformly defined over $s$.

\medskip
\noindent
{\bf Inside $\mathcal{B}_s\setminus \mathcal{T}_s$.} If $u\in \mathcal{B}_s\setminus \mathcal{T}_s$, then $\sup \tau _{s,u}<u$. Let $X_{s,u}^-=\mathcal{B}_{s,u}\cap (\sup \tau _{s,u},u)$ (recall we are assuming that $\mathcal{B}_{s,u}$  is an open subinterval with right endpoint $u$).  So, for all $v\in X_{s,u}^-$, since $\tau _{s,u}(v)\leq \sup \tau _{s,u}$, we have  $f_v(t)<f_u(t)$ for all $t\in (\sup \tau _{s,u}, v)$. From here we conclude that, for all $v\in X_{s,u}^-$ we have $C_{v|(\sup \tau _{s,u}, v)}\cap C_{u}=\emptyset $.  Hence (L) holds on $\mathcal{B}_s\setminus \mathcal{T}_s$.

\medskip
\noindent
{\bf Inside $\mathcal{T}_s$.} Without loss of generality, after removing  finitely many points we may assume that $\mathcal{T}_s$ is a finite union of open subintervals of $Y_s$ uniformly  defined over $s$.

Let $u\in \mathcal{T}_s$. Then

\begin{clm}\label{clm t_v increasing}
There is a maximal open subinterval of $\mathcal{B}_{s,u}$ with right endpoint $u$ on which $\tau _{s,u}$ is  strictly increasing and continuous.
\end{clm}

\pf
By the monotonicity theorem (\cite[Chapter 3, (1.2)]{vdd}), there is a maximal open subinterval of $\mathcal{B}_{s,u}$ with right endpoint $u$ on which $\tau _{s,u}$ is constant or strictly monotone and continuous. Now: 
\begin{itemize}
\item[-]
If $\tau _{u}$ is constant with value $t_0$, then $t_0=\tau _{s,u}(v)<v<u$ (Claim \ref{clm t_v<v}), which contradicts $\sup \tau _{s,u}=u$.

\item[-]
If $\tau _{u}$ is strictly decreasing, then for $v<w$ we have $\tau _{s,u}(w)<\tau _{s,u}(v)<v<u$ and the supremum of $\tau _{s,u}$ is the limit going to the left which leads again to a contradiction.
\end{itemize} 
\qed \\

By abuse of notation, we may assume that $\mathcal{B}_{s,u}$ is the maximal open subinterval with right endpoint $u$ on which $\tau _{s,u}$ is  strictly increasing and continuous. Since the family $\{\mathcal{B}_{s,u}\}_{u\in \mathcal{T}_s}$ is uniformly definable, we have a definable function $b: \mathcal{T}_s\to Y_s$, given by $b(u)=\inf \mathcal{B}_{s,u}$, such that $(b(u), u)=\mathcal{B}_{s,u}$. Consider the partition of $\mathcal{T}_s$ into finitely many points and open intervals given by Lemma \ref{lem mu}. Replace $\mathcal{T}_s$ by the collection of those finitely many open subintervals. It follows that given $u\in \mathcal{T}_s$, if  $w,v\in \mathcal{B}_{s,u}\cap  \mathcal{T}_s$ are such that $v<w$, then $v\in \mathcal{B}_{s,w}$. Moreover we have the following three possibilities:

\begin{itemize}
\item[-]
 $\tau _{s,w}(v)\leq \tau _{s,u}(v)<\tau _{s,u}(w)$;

\item[-]
$\tau _{s,u}(v)<\tau _{s,w}(v)<\tau _{s,u}(w)$;

\item[-]
$\tau _{s,u}(v)<\tau _{s,u}(w)\leq \tau _{s,w}(v)$.
\end{itemize}

\medskip
Suppose that there is an open subinterval $O$ of $\mathcal{T}_s$ defined over $s$ with no non-trivial points. Applying Lemma \ref{lem trivial and dim 2 family} to $O$, the definable function $b:O\to M$ and the definable function $\Omega :(b, \id _{O})_{O}\to M$ given by $\Omega (u,t)=\tau _{s,u}(t)$, we obtain that $O$ is a finite union of points and open subintervals, all defined over $s$,  such that on each such open subinterval $X$:
\begin{itemize}
\item[-]
there is a definable function $d_X :X\to M$ defined over $s$  such that  $b_{|X}(u)\leq d_X(u)<u$ for all $u\in X$;
\item[-]
for every  $u\in X$, if $v, w\in (d_X(u), u)\cap X$ are such that  $v<w$,  then we have $v\in (d_X(w), w)$ and $\Omega (u,v)=\Omega (w,v)$.
\end{itemize}

Given $u\in O$ let $X_{s,u}^-=(d_X(u), u)\cap X$ where $X$ as above is such that $u\in X$. Then for all $v,w\in X_{s,u}^-$, if $v<w$, then  $\tau _{s,w}(v)\leq \tau _{s,u}(v)<\tau _{s,u}(w)$. Now let $C_v$ (resp. $C_w$) denote the restriction $C_{v|(\tau _{s,u}(v),v)}$ (resp. $C_{w|(\tau _{s,u}(w),w)}$).  Since $f_{v|(\tau _{s,w}(v), v]}<f_{w|(\tau _{s,w}(v),v]}$, it follows that we have $C_v\cap C_w=\emptyset $ for all $v,w\in X_{s,u}^-$ with $v<w$.  Therefore $\{C_v\}_{v\in X_{s,u}^-}$ is a one-dimensional definable normal family of curves. Let $u'\in X$ be such that $u'<u$ and consider the restriction of the family $\{C_v\}_{v\in X_{s,u}^-}$ to $[u', u]\times [u', u]$. By \cite[Lemma 2.3]{pst},  the restriction of  $\{C_v\}_{v\in X_{s,u}^-}$ to $[u', u]\times [u', u]$ is definable over $[u',u]$, hence over $O$. So by Theorem \ref{thm Z1 and curves} applied to ${\mathbb M}_{|O}$ there is a non-trivial point in $O$, which is a contradiction, since we assumed that $O$ had no non-trivial points.
\qed \\

In the next lemma, although the condition (N2) is exactly the same as the condition (N1) from Lemma \ref{lem type (N1) or type (L)}, the curves appearing in (N2) are constructed in a different way than those appearing in (N1) and it is convenient for now to keep track of this information. A similar remark applies to condition (D1) from Lemma \ref{lem type (N1) or type (L)} and (D2) below.

\begin{lem}\label{lem type (L)} 
Let $X\in \mathcal{X}_s$ be  of type (L). Then after removing finitely many points from $X$, all uniformly defined over $s$, and replacing $X$ by one of the remaining open subintervals we may assume that the following holds
\begin{itemize}
\item[(N2)] 
for every $u\in X$  there is an open subinterval $X_{s,u}^-\subseteq Y_{s,u}^-$  with right endpoint $u$ such that  for every  $v\in X_{s,u}^-$ there is an open subinterval $J_{s,u,v}^-\subseteq I_{s,v}^-$ with right endpoint $v$ such that $\{C_{v|J_{s,u,v}^-}\}_{v\in X_{s,u}^-}$ is a one-dimensional definable  normal family of  curves,
\end{itemize}
or else
\begin{itemize}
\item[(D2)] 
every open subinterval of $X$ defined over $s$ has a non-trivial point.
\end{itemize}
\end{lem}

\pf
Let $X\in \mathcal{X}_s$ be of type (L). Since the family $\{X_{s,u}^-\}_{u\in X}$ of open subintervals of $X$, each with right end point $u$,  is uniformly definable, there is  a definable function $\mu :X \to Y_s$, given by $\mu (u) =\inf X_{s,u}^-$, such that  $(\mu (u),u)=X_{s,u}^-$.  Now consider the partition  of $X$  into  finitely points and open intervals given by Lemma \ref{lem mu}. Note that, since the family $\{\mu :X\to Y_s\}$ is uniformly defined over $s$ and the conclusion of Lemma \ref{lem mu} is first-order in $s$, these finitely many points and open intervals are uniformly defined over $s$. Replace $X$ by one such open subinterval. It follows that if $u\in X$, then for all $v,w\in X_{s,u}^-$ with  $v<w$ we have $v\in X_{s,w}^-$.

\medskip
\noindent
Fix $u\in X$. If  $v,w\in X_{s,u}^-$ are such that $v<w$, then, since $v\in X_{s,w}^-$, there are open subintervals  $J_{s,w, v}^-$ of $I_{s,v}^-$ and $J_{s,w}^-$ of $I_{s,w}^-$, uniformly in $v$ and $w$,  with right endpoint $v$ and $w$ respectively such that  $C_{v|J_{s,w,v}^-}\cap C_{w|J_{s,w}^-}$ is finite, equivalently $\dim (C_{v|J_{s,w, v}^-}\cap C_{w|J_{s,w}^-})=0$. So there is  a uniformly definable family of definable functions 
$$\{j_{s,v}:(v,u)\to I_{s,v}^-\}_{v\in X_{s,u}^-}$$ 
given by  $j_{s,v}(w)=\inf J_{s,w,v}^-$ for $w\in (v,u)$. 

By construction we have  $J_{s,w,v}^-=(j_{s,v}(w),v)$, in particular $j_{s,v}(w)<v$, for all $w\in (v,u)$ and all $v\in X_{s,u}^-$. 

Let
$$\mathcal{V}_{s,u}=\{v\in X_{s,u}^-: \sup j_{s,v}=v\}.$$

Then for each $u\in X$ either $\mathcal{V}_{s,u}$ or $X_{s,u}^-\setminus \mathcal{V}_{s,u}$, but not both, contains an open subinterval of $X_{s,u}^-$ with right endpoint $u$. 
Let 

$$\mathcal{V}_s=\{u\in X: \mathcal{V}_{s,u}\,\,\textrm{contains an open subinterval with right endpoint}\,\, u\}.$$

Note that $\mathcal{V}_s$ is uniformly defined over $s$.

\medskip
\noindent
{\bf Inside $X\setminus \mathcal{V}_s$.} We have that $X\setminus \mathcal{V}_s \subseteq Y_{s}$ is a finite union of points and open subintervals of $Y_s$, all uniformly defined over $s$. Removing these finitely many points we may assume that $X\setminus \mathcal{V}_s$ is a finite union of open subintervals of $Y_s$ uniformly defined over $s$. 

If $u\in X\setminus \mathcal{V}_s$, then $X_{s,u}^-\setminus \mathcal{V}_{s,u}$ contains an open subinterval with right endpoint $u$. Let us assume, with abuse of notation, that $X_{s,u}^-$ is the maximal  open such subinterval. If $v\in X_{s,u}^-$, then $\sup j_{s,v}<v$.  Note also that $v\in X$. Consider the open subinterval $(\sup j_{s,v}, v)\cap J^-_{s,v}$, which by abuse of notation we still denote by $J_{s,u,v}^-$. Then for all $w\in (v,u)$  we have $\dim (C_{v|J_{s,u,v}^-}\cap C_{w|J_{s,w}^-})=0$ since $j_{s,v}(w)<\sup j_{s,v}$. Now if $w\in X_{s,u}^-$ is such that $v<w$, then since $(\sup j_{s,w}, w)\cap J^-_{s,w}\subseteq J^-_{s,w}$ we have  $\dim (C_{v|J_{s,u,v}^-}\cap C_{w|J_{s,u,w}^-})=0$. Hence $\{C_{v|J_{s,u,v}^-}\}_{v\in X_{s,u}^-}$ is a one-dimensional definable  normal family of  curves, and  (N2) holds on $X\setminus \mathcal{V}_s$.

\medskip
\noindent
{\bf Step 2: Inside $\mathcal{V}_s$.} If $u\in \mathcal{V}_s$, then $\mathcal{V}_{s,u}$ contains an open subinterval with right endpoint $u$. Let us assume, with abuse of notation, that $\mathcal{V}_{s,u}$ is the maximal  open such subinterval. 

Let $v\in \mathcal{V}_{s,u}$. Then  the following holds.

\begin{clm}\label{clm D_v for L}
 There is a maximal open subinterval $D_{s,v}$ of $(v,u)$ such that $j_{s,v|D_{s,v}}:D_{s,v}\to I_{s,v}^-$ is strictly decreasing and $D_{s,v}$ has left endpoint $v$.
\end{clm}
\pf 
By the monotonicity theorem (\cite[Chapter 3, (1.2)]{vdd}), $(v,u)$ is a finite union of points  and open subintervals such that on each such open subinterval the restriction of $j _{s,v}$ to that interval is either constant or strictly monotone and continuous.  Since $\sup j_{s,v}=v$, there is an open subinterval $D_{s,v}$ of $(v,u)$ such that $j_{s,v|D_{s,v}}:D_{s,v}\to I_{s,v}^-$ is strictly monotone and continuous and $j_{s,v}(D_{s,v})$ is an  open subinterval of $I_{s,v}^-$ with left endpoint $v$.

Let $a^-$ (resp. $a^+$) be the left (resp. right) endpoint of $D_{s,v}$. Note  that there are open subintervals of $D_{s,v}$ with right endpoint $a^+$ (resp. left endpoint $a^-$)  such that the restriction of the definable function $a\mapsto H_s(a,j_{s,v}(a))$ of Lemma \ref{lem H_s} is continuous on them.  Note also that for all $a\in D_{s,v}$ we have $f_{v}(j_{s,v}(a))=f_a(j_{s,v}(a))$ since, by definition of $j_{s,v}(a)$,  the curves $C_v$ and $C_a$ coincide in a subinterval with right endpoint $j_{s,v}(a)$. 

Suppose that  $j_{s,v|D_{s,v}}$ is strictly increasing, in particular, $\lim _{a\to a^+}j_{s,v}(a)=\sup j_{s,v}=v$. Then by continuity, we have 
\begin{eqnarray*}
v & = & \lim _{a\to a^+}f_v(j_{s,v}(a))\\
   &= & \lim _{a\to a^+}f_a(j_{s,v}(a))\\
   & = & \lim _{a\to a^+}H_s(a, j_{s,v}(a))\\
   & =& f_{a^+}(v)
\end{eqnarray*}
contradicting the fact that $v<f_{a^+}(v)$ since $v<a^+$.

Suppose that  $j_{s,v|D_{s,v}}$ is strictly decreasing but $v<a^-$. Then by continuity again, we have $v=\lim _{a\to a^-}f_v(j_{s,v}(a))=\lim _{a\to a^-}f_a(j_{s,v}(a))=f_{a^-}(v)$ contradicting the fact that $v<f_{a^-}(v)$ since $v<a^-$. \qed \\

Since the family $\{D_{s,v}\}_{v\in \mathcal{V}_{s,u}}$ is uniformly definable, we have a definable function $b: \mathcal{V}_{s,u}\to Y_s$, given by $b(v)=\sup D_{s,v}$, such that $(v, b(v))=D_{s,v}$. Consider the partition of $\mathcal{V}_{s,u}$ into finitely many points and open intervals given by Lemma \ref{lem mu+}. Replace $\mathcal{V}_{s,u}$ by the maximal rightmost such open subinterval. It follows that given $v\in \mathcal{V}_{s,u}$, if  $x,y\in D_{s,v}\cap  \mathcal{V}_{s,u}$ are such that $x<y$, then $x\in D_{s,y}$. Moreover we have the following three possibilities:

\begin{itemize}
\item[-]
$j_{s,x}(y)\leq j_{s,v}(y)< j_{s,v}(x)$;

\item[-]
$j_{s,v}(y)<j_{s,x}(y)<j_{s,v}(x)$;

\item[-]
$j_{s,v}(y)<j_{s,v}(x)\leq j_{s,x}(y)$.
\end{itemize}

\medskip
Suppose that there is an open subinterval $O$ of $\mathcal{V}_s$ defined over $s$ with no non-trivial points. Then for every $u\in O$, $\mathcal{V}_{s,u}$ has no non-trivial points. Applying Lemma \ref{lem trivial and dim 2 family+} to $\mathcal{V}_{s,u}$, the definable function $b:\mathcal{V}_{s,u}\to M$ and the definable function $\Omega :(\id _{\mathcal{V}_{s,u}}, b)_{\mathcal{V}_{s,u}}\to M$ given by $\Omega (v,t)=j_{s,v}(t)$, we obtain that $\mathcal{V}_{s,u}$ is a finite union of points and open subintervals, all defined over $s$,  such that on each such open subinterval $X$:
\begin{itemize}
\item[-]
there is a definable function $d_X :X\to M$ defined over $s$  such that  $v<d_{X}(v)\leq b_{|X}(v)$ for all $v\in X$;
\item[-]
for every  $v\in X$, if $x, y\in (v, d_X(v))\cap X$ are such that  $x<y$,  then we have $y\in (x, d_X(x))$ and $\Omega (v,y)=\Omega (x,y)$.
\end{itemize}

Let $X_{s,u}^-$ be the open subinterval $X$ of  $\mathcal{V}_u$  with right endpoint $u$. Given $v\in X_{s,u}^-$ let $X_{s,u,v}^+=(v, d_X(v))\cap \mathcal{V}_u$. Then for all $x, y\in X_{s,u,v}^+$, if $x<y$, then  $j_{s,x}(y)\leq j_{s,v}(y)< j_{s,v}(x)$. Now let $C_x$ (resp. $C_y$) denote the restriction $C_{x|(j_{s,v}(x),x)}$ (resp. $C_{y|(j_{s,v}(y),y)}$).  Since $C_{x|(j_{s,x}(y), x)}\cap C_{y|(j_{s,x}(y),x)}$ is finite, it follows that $C_x\cap C_y$ is finite for all $x, y\in X_{s,u,v}^+$ with $x<y$.  Therefore $\{C_x\}_{x\in X_{s,u,v}^+}$ is a one-dimensional definable normal family of curves. Let $v'\in X$ be such that $v<v'$ and consider the restriction of the family $\{C_x\}_{x\in X_{s,u,v}^+}$ to $[v, v']\times [v, v']$. By \cite[Lemma 2.3]{pst},  the restriction of  $\{C_x\}_{x\in X_{s,u,v}^+}$ to $[v, v']\times [v, v']$ is definable over $[v,v']$, hence over $O$. So by Theorem \ref{thm Z1 and curves} applied to ${\mathbb M}_{|O}$ there is a non-trivial point in $O$, which is a contradiction since we assumed that $O$ had no non-trivial points.  \qed \\

The cost for obtaining Lemmas \ref{lem type (N1) or type (L)} and \ref{lem type (L)} was the introduction of extra parameters. To proceed we will apply Lemma \ref{lem removing parameters} to remove first the  parameter $s$ (Lemma \ref{lem uniform type (N) or else}) and then the extra parameters (Lemma \ref{lem main 0}).

\begin{lem}\label{lem uniform type (N) or else} 
Consider the $\emptyset$-definable open intervals  $I\in \mathcal{I}$ and $Y\in \mathcal{Y}^I$ of Lemmas \ref{lem def skolem and family of funct} and \ref{lem H_s}. Then after naming finitely many elements of $\mathbb{P}$, and after removing finitely many $\emptyset $-definable points and replacing each $Y$ by one of the remaining open subintervals we may assume that for each $Y\in \mathcal{Y}^I$ there are $\emptyset $-definable functions $q^-, q^+:Y\to M$  with $q^-<q^+$, such that  $y\in (q^-(y), q^+(y))$ for all  $y\in Y$. Moreover, letting $U_y=(q^-(y), q^+(y))$, then for all $y\in Y$ the following holds
\begin{itemize}
\item[(N$'$)] 
for every $u\in U_y$  there is an open subinterval $X_{y,u}^-\subseteq U_y$ with right endpoint $u$, uniformly definable in $y$,  such that  for every  $v\in X_{y,u}^-$ there is an open subinterval $J_{y,u,v}^-\subseteq U_y$ with right endpoint $v$, uniformly definable in $y$,  such that $\{C_{v|J_{y,u,v}^-}\}_{v\in X_{y,u}^-}$ is a one-dimensional definable  normal family of  curves,
\end{itemize}
or else for all $y\in Y$ 
\begin{itemize}
\item[(D$'$)]
every open subinterval of $U_y$ defined over $y$ has a non-trivial point.
\end{itemize}
\end{lem}

\pf
Consider the $\emptyset$-definable open intervals  $I\in \mathcal{I}$, $Y\in \mathcal{Y}^I$ and  the definable intervals
\begin{eqnarray*}
Y_s=
\begin{cases}
Y\cap (s,e_Y)\qquad \textrm{if}\,\,\, s<e_Y\\
\\
Y\cap (e_Y,s)\qquad \textrm{otherwise}
\end{cases}
\end{eqnarray*}
 with $s\in S$, introduced in Lemmas \ref{lem def skolem and family of funct} and \ref{lem H_s}. Consider also the  finite collection $\mathcal{X}_s$,  uniformly defined over $s$, of two by two disjoint  open intervals  defined over $s$, given by Lemmas \ref{lem type (N1) or type (L)} and \ref{lem type (L)}, such that $Y_s\setminus \bigcup \mathcal{X}_s$ is finite and for each $X\in \mathcal{X}_s$ either $X$ is of one of the types (N1), (D1), (N2) or (D2). 

If consider the data $(Y, S, \{Y_s\}_{s\in S}, \{\mathcal{X}_s\}_{s\in S})$ then we are in the set up of Lemma \ref{lem removing parameters}. So after naming finitely many elements of $\mathbb{P}$, and after removing  finitely many $\emptyset $-definable points from $Y$ and replacing $Y$ by one of the remaining open subintervals, we may assume that there are $\emptyset$-definable functions $q^-, q^+:Y\to M$  such that $q^-<q^+$ and $s:Y\to S$ such that, for each $y\in Y$,  $(q^-(y), q^+(y))$ is one of the open subintervals of $Y_{s(y)}$ in $\mathcal{X}_{s(y)}^-$ and $y\in (q^-(y), q^+(y))$.

By the remarks about uniform definability over $s(y)$ throughout the proof of Lemma \ref{lem type (N1) or type (L)}, we see that the open subinterval $(q^-(y), q^+(y))$ being of type (N1), not being of type (N1) but being of type (L) and not being of type (N1) neither of type (L) but of type (D1) are first-order properties of  $y$. Indeed, these correspond to: the open subinterval $(q^-(y), q^+(y))$ is contained in $Y_{s(y)}\setminus \mathcal{B}_{s(y)}$ (resp. in $\mathcal{B}_{s(y)}\setminus \mathcal{T}_{s(y)}$ or in $\mathcal{T}_{s(y)}$). So one can partition $Y$ into one, two or three $\emptyset $-definable subsets  such that  on each one of them,  all the open subintervals $(q^-(y), q^+(y))$ are of type (N1), or all are of type (L) but not of type (N1) or none is of type (N1) or type (L) but are of type (D1). 

Similarly, using the remarks about uniform definability over $s(y)$ throughout the proof of Lemma \ref{lem type (L)}, we see that the part of $Y$ corresponding to type (L) but not type (N1) can be further partitioned into one or two  $\emptyset $-definable subsets  such that  on each one of them, all the open subintervals $(q^-(y), q^+(y))$ are of type (N2)  or none is of type (N2) and all are of type (D2). 

So removing finitely many $\emptyset $-definable points and replacing $Y$ by one of the remaining open subintervals, we can assume that the open intervals $(q^-(y), q^+(y))$ are either all  of type (N1), or all of type (D1),  or all of type (N2) or all of type (D2). 

Using the fact that for each $y\in Y$,  $(q^-(y), q^+(y))$ is an open subinterval of $Y_{s(y)}$,  $I_{s(y), v}^-$ is an open subinterval with right endpoint $v$ and $v\in (q^-(v), q^+(v))$ for each $v\in (q^-(y), q^+(y))$, the result now follows at once from Lemmas \ref{lem type (N1) or type (L)} and \ref{lem type (L)}. \qed \\

We can  apply Lemma \ref{lem removing parameters} once more to remove the  parameter $u$ from  Lemma \ref{lem uniform type (N) or else} and reach our main technical lemma:

\begin{lem}\label{lem main 0}
After naming finitely many elements of $\mathbb{P}$, there is a collection $\mathcal{Y}$ of finitely many $\emptyset $-definable open subinterval of $M$ such that $M\setminus \bigcup \mathcal{Y}$ is finite and for each $Y\in \mathcal{Y}$ there is a uniformly  $\emptyset $-definable family $\{X_y\}_{y\in Y}$ of open subintervals such that $y\in X_y$ for all  $y\in Y$ and moreover, for all $y\in Y$ the following holds
\begin{itemize}
\item[(N)] 
for every  $v\in X_{y}$ there is an open subinterval $J_{y,v}^-$ of $X_y$ with right endpoint $v$, uniformly definable in $y$,  such that $\{C_{v|J_{y,v}^-}\}_{v\in X_{y}}$ is a one-dimensional definable  normal family of  curves,\\
\end{itemize}
or else for all $y\in Y$
\begin{itemize}
\item[(D)]
every open subinterval of $X_y$ defined over $y$ has a non-trivial point.
\end{itemize}
\end{lem}

\pf
We start with $I\in \mathcal{I}$,  $Y\in \mathcal{Y}^I$ and  $U_y$ as in Lemma \ref{lem uniform type (N) or else}.  

\medskip
\noindent
Case (D$'$): Suppose that   for all $y\in Y$, $U_y$ is of type (D$'$). Then  set $X_y=U_y$. \\

\medskip
\noindent
Case (N$'$): Suppose that for all $y\in Y$, $U_y$ is of type (N$'$). Let $Z=Y$, $T=Y$ and for $t\in T$ let $Z_t=U_t$. Note that $F_t:=Z_t \setminus \bigcup \{X_{t,u}^-:u\in Z_t\}$ is finite. Otherwise, this definable subset of $M$ would contain an open subinterval $W$. But then, for  $u\in W$, we would have $W\cap X_{t,u}^-\neq \emptyset $ since $W\subseteq X$ is an open subinterval and, on the other hand,  $W\cap X_{t,u}^-=\emptyset $ since $W\subseteq F_t$. 

Let $\mathcal{U}_t$ be the collection of two by two disjoint open subintervals obtained by removing $F_t$ from $Z_t$.  If we consider the data $(Z,T, \{Z_t\}_{t\in T}, \{\mathcal{U}_t\}_{t\in T})$, then we are in the set up of Lemma \ref{lem removing parameters}. Therefore,  after naming finitely many elements of $\mathbb{P}$, and after removing  finitely many $\emptyset $-definable points from $Z$ and replacing $Z$ by one of the remaining open subintervals, we may assume that there are $\emptyset$-definable functions $p^-, p^+:Z\to M$  such that $p^-<p^+$ and $t:Z\to T$ such that, for each $z\in Z$,  $(p^-(z), p^+(z))$ is one of the open subintervals of $Z_{t(z)}$ in $\mathcal{U}_{t(z)}$ and $z\in (p^-(z), p^+(z))$. Hence, for each $z\in Z$, there is $u\in Z_{t(z)}$ such that $z\in X_{t(z), u}^-$. By definable Skolem functions there is another definable function $t':Z\to T$ such that, for each $z\in Z$, we have $t'(z)\in Z_{t(z)}\subseteq T$ and $z\in X_{t(z), t'(z)}^-\subseteq Z_{t(z)}\subseteq Z$. Varying the parameter over which $t':Z\to T$ is defined, we get a uniformly $\emptyset $-definable family of such functions. By Proposition \ref{prop def skolem dcl} we can choose a parameter in ${\rm dcl}(p')$ for some $p'\in P^k$. After naming also the elements in $p'$ we see that we assume that $t': Z\to T$  is defined over $\emptyset $.

For $z\in Z$, take $X_z=X_{t(z), t'(z)}^-$ and for $v\in X_z$ take $J_{z,v}^-=J_{t(z), t'(z), v}^-$. \\

It is now clear that if we let $\mathcal{Y}$ be the collection of all those finitely many $Y$'s in $\mathcal{Y}^I$ from Case (D$'$) together with those finitely many $Z$'s obtained in Case (N$'$), the result follows. 
\qed \\

\end{section}

\begin{section}{The main results}\label{section main results}
{\it Here if we don't say otherwise,  $\mathbb{M}$ is a sufficiently saturated o-minimal structure with definable Skolem functions.}
\medskip

The first consequence of the main technical lemma (Lemma \ref{lem main 0}) is the density of non-trivial points in $M$:

\begin{prop}\label{prop non-trivial are dense} 
Every open subinterval of $M$ defined over $s$ has, after naming finitely many elements of $\mathbb{P}$,  a non-trivial point over $s$. In particular, the set of non-trivial points is dense in $M$.
\end{prop}

\pf
Let $O$ be an open subinterval defined over $s$. Note if we add the parameters $s$ to the language, the structure ${\mathbb M}$ in the new language still has definable Skolem functions. So it is enough to proceed  assuming that $O$ is defined over $\emptyset $. 

Let $\mathcal{Y}$  be as in Lemma \ref{lem main 0}. Then it is enough to show that for each $Y\in \mathcal{Y}$, $O\cap Y$ has a non-trivial point over $\emptyset $. Fix $Y\in \mathcal{Y}$. By Proposition \ref{prop def skolem dcl}, after naming finitely many elements of $\mathbb{P}$, we can choose  $y\in O\cap Y$ defined over $\emptyset $.   

\medskip
\noindent
Case (N): Suppose that $Y$ is of type (N). Take $u'<u$ both in $O\cap Y$ such that $y\in (u',u)$.  By \cite[Lemma 2.3]{pst},  the restriction of  $\{C_{v|J^-_{y,v}}\}_{v\in X_{y}}$ to $[u', u]\times [u', u]$ is definable over $[u',u]$, hence over $O\cap Y$. So by Theorem \ref{thm Z1 and curves} applied to ${\mathbb M}_{|O\cap Y}$ there is a non-trivial point in $O\cap Y$.  

\medskip
\noindent
Case (D): Suppose that $Y$ is of type (D). Then $O\cap Y\cap X_y$ is an open sub-interval of $X_y$ defined over $y$. So there is a non-trivial point in $O\cap Y$. 

\medskip
Now we conclude using:

\begin{clm}\label{clm non trivial point over s}
If $O$ is an open subinterval defined over $s$ with a non-trivial point, then after naming finitely many elements of $\mathbb{P}$ it has a non-trivial point over $s$.
\end{clm}

\pf
Let $a$ be a non-trivial point in $O$. Since $O$ is an open interval, by \cite[Lemma 2.5]{pst}, $a$ is non-trivial over $O$. Therefore, there is an ${\mathbb M}_{|O}$-definable map $G:I\times I\to J$ with $I, J\subseteq O$ open subintervals, which is continuous and strictly monotone in each variable. Let $c\in O^m$ be the parameter over which $G$ is defined. Then there is a uniformly definable family $\{G_u:I_u\times I_u\to J_u\}_{u\in U}$, defined over $s$, of definable maps such that for each $u\in U\subseteq O^m$, $I_u, J_u\subseteq O$ are open subintervals, $G_u:I_u\times I_u\to J_u$ is continuous and strictly monotone in each variable,  $c\in U$ and $G_c=G$. By Proposition \ref{prop def skolem dcl} and after naming finitely many elements of $\mathbb{P}$ we can  choose $u\in U$ defined over $s$ and  then choose $b\in I_u$ defined over $s$. This $b$ is therefore a non-trivial point of $O$ over $s$. \qed 

\qed \\

Density of non-trivial points in $M$ implies the following:

\begin{lem}\label{lem trivial are over emptyset}
Every trivial point of $M$ is $\mathbb{P}$-definable.
\end{lem}

\pf
Let $t$ be a trivial point of $M$ and assume that it is defined over $c\,\mathbb{P}$ for some $c=(c_1, \ldots, c_l)\in M^l$. Then there is a $\mathbb{P}$-definable function $f:Z\subseteq M^l\to M$ such that $c\in Z$ and $f(c)=t$. We can assume that there is no $c'\subset c$ for which there is  a $\mathbb{P}$-definable $f':Z'\subseteq M^{k}\to M$ such that $c'\subseteq Z'$ and $f'(c')=t$. Note that if some $c_i$ is defined over $c'\,\mathbb{P}$ for  $c'\subset c$ consisting of the remaining elements, then we may replace $f:Z\to M$ by another $\mathbb{P}$-definable $f':Z'\subseteq M^{l-1}\to M$ such that $c'\subseteq Z'$ and $f'(c')=t$. So no $c_i$ is defined over $c'\,\mathbb{P}$ for $c'\subset c$ consisting of the remaining elements. 

Since ${\rm dim}(c/\mathbb{P})=\dim Z$, by the cell decomposition theorem, we can furthermore assume that $Z$ is an open cell in $M^l$. 

Let $\pi :M^l\to M^{l-1}$ be the projection onto the first $l-1$ coordinates. Let $W=\{w' \in \pi (Z): (w',c_l)\in Z\}$. Since $Z$ is open, $W$ is an open definable subset of $\pi (Z)$, defined over $c_l\, \mathbb{P}$. For each $w'\in W$, let $I_{w'}=\{x\in M: (w',x)\in Z\}$ and let $g_{w'}:I_{w'}\to M$ be the definable function given by $g_{w'}(x)=f(w',x)$. Then each $I_{w'}$ is an open subinterval of $M$ with $c_l\in I_{w'}$ and $g_{w'}:I_{w'}\to M$ is uniformly defined over $w' \, \mathbb{P}$.  

Let $c'=(c_1, \ldots, c_{l-1})\in W$. Then by the monotonicity theorem, $I_{c'}$ is a finite union of points and open subintervals all defined over $c'\, \mathbb{P}$ such that on each such open subinterval $g_{c'}$ is either constant or strictly monotone and continuous. Note that $c_l$ cannot be one the these points, since otherwise it would be defined over $c'\,\mathbb{P}$. So $c_l$ is in one of the open subintervals, say $J$. Then $g_{c'|J}$ must be strictly monotone and continuous, otherwise it would be constant and, by  Proposition \ref{prop def skolem dcl},  we could choose  $b\in J$ defined over $c'\,\mathbb{P}$ and obtain another $\mathbb{P}$-definable function $f':Z'\subseteq M^{l-1}\to M$  such that $c'\in Z'$ and $f'(c')=t$, contradicting our assumption on $c$ and $t$. 

Now the properties $J\subseteq I_{c'}$ is an open subinterval, $c_l\in J$ and $g_{c'|J}$ is strictly monotone and continuous are first-order in $c'$. So there is a definable subset $U\subseteq W$, defined over $c_l\,\mathbb{P}$, and a uniformly definable family of open subintervals $J_{u'}\subseteq I_{u'}$ such that: (i) $c'\in U$; (ii) $J_{c'}=J$; (iii) for each $u'\in U$ we have $c_l\in J_{u'}$ and $g_{u'|J_{u'}}$ is strictly monotone and continuous. Moreover, since $c'\in U\subseteq W$ and ${\rm dim}(c'/\mathbb{P} )=\dim W$, by the cell decomposition theorem we may assume that $U$ is an open cell in $M^{l-1}$.

Consider the definable function $U\to M:u'\mapsto g_{u'|J_{u'}}(c_l)$, defined over $c_l\,\mathbb{P}$. Its image, which is defined over $c_l\,\mathbb{P}$,  must be infinite, otherwise $t$, which is in the image, would be defined over $c_l\,\mathbb{P}$. Take one of the maximal open subintervals contained in this image, say $K$. Since $g_{c'|J}$ is strictly monotone and $g_{c'|J}(c_l)=t$, $c_l$ is a trivial point of $M$. So similarly, every point in $K$ must be a trivial point, contradicting Proposition \ref{prop non-trivial are dense}. 
\qed \\

We are ready to show the main result:

\begin{thm}\label{thm main 0}
Let $\mathbb{M}$ be an o-minimal structure. Then the following are equivalent:
\begin{enumerate}
\item[(1)]
$\mathbb{M}$ has definable Skolem functions. 
\item[(2)]
After naming finitely many elements of $\mathbb{P}$, there is a finite collection $\mathcal{Y}$ of $\emptyset $-definable open subintervals of $M$ such that:
\begin{itemize}
\item[-]
$M\setminus \bigcup \mathcal{Y}$ is a finite set of trivial points each defined over $\emptyset $.
\item[-]
each $Y\in \mathcal{Y}$ is the union of a uniformly $\emptyset $-definable family of group-intervals, each with a fixed positive element, parametrized by the end points of the intervals.
\item[-]
each $Y\in \mathcal{Y}$ has  a fixed $\emptyset $-definable element.
\end{itemize}
\item[(3)] $\mathbb{M}$  has definable choice.
\end{enumerate}
\end{thm}

\pf
First note that if the result holds in a sufficiently saturated elementary extension of $\mathbb{M}$, then it holds for $\mathbb{M}$. So we may assume that $\mathbb{M}$ sufficiently saturated.

(1) $\Rightarrow $ (2): 
If $x\in M$ is non-trivial, then by the Trichotomy theorem $x$ is contained in a definable group-interval  $(I, 0, +, <)$. As observed in \cite[Note 3.2]{epr}, by varying the parameter over which $+$ and $0$ are defined but requiring that the domain remains $I$, there is a uniformly definable family of such $+$'s and corresponding $0$'s defined over the end points of $I$. By Proposition \ref{prop def skolem dcl}, we can choose a parameter  in the family defined over the end points of $I$ and some $p\in P^k$. Therefore, after naming the elements in $p$, we may assume that $+$ and $0$ are defined over those end points. Applying again Proposition \ref{prop def skolem dcl} and naming finitely many the elements of $\mathbb{P}$ we can choose a positive element $1$ in $(I,0,+,<)$ defined over the end points of $I$. 

Suppose that $I=(c_1, c_2)$ and let $c=(c_1, c_2)$.  Varying $c$ we obtain a $\emptyset $-definable family $\{ (I_w, 0_w, 1_w, +_w,<)\}_{w\in W}$ of definable group-intervals, each with a fixed positive element, such that $x\in I_W$ where $I_W=\bigcup _{w\in W}I_w$. 

Let $\mathcal{W}$ be the set of all such $I_W$'s and let $\mathcal{T}$ be the set of all trivial points of $M$. Then $\mathcal{T}\cup \bigcup \mathcal{W}=M$. Since  every element of $\mathcal{T}$ is defined over $\mathbb{P}$ (Lemma \ref{lem trivial are over emptyset}),  by saturation, there is a finite subset $T\subseteq \mathcal{T}$ and finitely many $I_{W_1}, \ldots , I_{W_l}\in \mathcal{W}$ such that $T\cup I_{W_1}\cup  \ldots \cup  I_{W_l}=M$. Since $\mathcal{T}\cap (I_{W_1}\cup  \ldots \cup  I_{W_l})=\emptyset $ and,  by Proposition \ref{prop def skolem dcl}, for each $i$, after naming finitely many the elements of $\mathbb{P}$, we can choose a $\emptyset $-definable element in the $\emptyset $-definable set $I_W$, the result follows.

(2)  $\Rightarrow $ (3): Suppose that (2) holds. Then we have $T\sqcup I_{W_1}\sqcup  \ldots \sqcup  I_{W_l}=M$ where $T$ is a finite set of trivial points of $M$ all defined over $\mathbb{P}$ and each $I_W$ is of the form $I_W=\bigcup _{w\in W}I_w$ with $\{ (I_w, 0_w, 1_w, +_w,<)\}_{w\in W}$ a $\mathbb{P}$-definable family of definable group-intervals, each with a fixed positive element, with $w$ the pair of left and right end point of $I_w$ in $M\cup \{\pm \infty \}$. Moreover, for each $i$, let $e_i\in I_{W_i}$ be a fixed $\mathbb{P}$-definable element. 

Let   $e=(e_1, \ldots, e_l)$. Let $X\subseteq M$ be a non empty definable subset of $M$ defined over $a=(a_1,\ldots, a_n)$. By Proposition \ref{prop def choice dcl} we need to show that there is an element of $X$ in ${\rm dcl}([a]e\,\mathbb{P})={\rm dcl}([a]\,\mathbb{P})$.  

Take a cell decomposition of $M$ respecting $X$ and the $\mathbb{P}$-definable finite set $X\cap T$. Then $X$ is the disjoint  union of points $d_1< \ldots < d_k$  containing $X\cap T$ and open intervals $J_1<\ldots <J_m$  with the $d_i$'s and the end points of the $J_j$'s in ${\rm dcl}(a\,\mathbb{P})$.  By varying $a$ we have a uniformly $\mathbb{P}$-definable family $\{X_s\}_{s\in S}$ of subsets of $M$ such that $a\in S$ and $X_a=X$, together with a finite collection $\mathcal{G}$ of  $\mathbb{P}$-definable functions from $S$ to $M$ such that: (i) for each $s\in S$, $X_s$ is the disjoint union of $k$ many points and $m$ many open intervals; (ii) for each $g\in \mathcal{G}$ and $s\in S$, $g(s)$ is one of the $k$ many points of $X_s$ exactly when $g(a)$ is one of the $k$ many points of $X$ or $g(s)$ is one of end points of the $m$ many of open subintervals of $X_s$ exactly when $g(a)$ is an end point of one of the $J_j$'s; (iii) for every $g, g'\in \mathcal{G}$ and $s\in S$ we have $g(s)<g'(s)$ if and only if $g(a)<g'(a)$. 

Note that if $s\in [a]$, then  $X_s=X$ and so for each $g\in \mathcal{G}$ we have $g(s)=g(a)$. Therefore, the $g$'s show that each of  the $d_i$'s and each of the end points of the open subintervals $J_j$'s are in ${\rm dcl}([a]\,\mathbb{P})$.

If $k\geq 1$ then we are done. So we may assume that $X$ is disjoint union of open intervals $J_1<\ldots <J_m$ and $X\cap T=\emptyset $. Let $J=J_1$ and take $W\in \{W_1, \ldots , W_l\}$ the unique element such that $J\subseteq I_W$. Then $J$ has end points in $\cl (I_W)$ and either $J$ has an end point in $I_W$ or $J=I_W$. If $J=I_W$, then one of the elements of $\{e_1, \ldots, e_l\}$  is in  $X\cap {\rm dcl}([a]e\,\mathbb{P})$. 

Suppose $J$ has an end point $c$ in $I_W$ (which is  in ${\rm dcl}([a]\,\mathbb{P})$). Let $g\in \mathcal{G}$ be the unique element such that $g(a)=c$. Let $U=\{u\in W: c\in I_u\}$,   

\[
v_1=\inf \{x:x<c\,\,\textrm{and}\,\, (x, y)\in U\,\,\textrm{for some}\,\,c<y\},
\]
\[
v_2=\sup \{y: c<y,\textrm{and}\,\, (v_1, y)\in U\}
\]
and  $v=(v_1,v_2)$. Then $U$ is definable over $c\,\mathbb{P}$, $v\in U$ and is also defined over $c\,\mathbb{P}$. So there are $\mathbb{P}$-definable functions $h_1, h_2:C\to M$ such that $c\in C$ and $h_i(c)=v_i$ for $i=1,2$. Then the $\mathbb{P}$-definable functions $h_{i|g(S)\cap C}\circ g: S\to M$ ($i=1,2$) show that $v\in {\rm dcl}([a]\,\mathbb{P})$.

Since $c\in I_v$, we have that $J\cap I_v$ is an open subinterval of $I_v$. Let $c'$ be the other end point of $J\cap I_v$ in $\cl (I_v)$. Then either $c'\in I_v$ and in this case $J\subseteq I_v$ or $c'$ is one of the end points of $I_v$. In both cases we have $c'\in {\rm dcl}([a]\,\mathbb{P})$. 

Suppose that $c'\in I_v$.  In $\mathbb{M}_{|I_v}$, the set $J\cap I_v$ is definable over $c\,c'$. Since $\mathbb{M}_{|I_v}$ is an o-minimal expansion of a group-interval with a fixed positive element, it has definable choice  (as in \cite[Fact 4.5]{epr} or \cite[Chapter 6 (1.1)]{vdd}). Therefore, by Proposition \ref{prop def choice dcl} in the o-minimal structure $\mathbb{M}_{|I_v}$, there is an element of $J$ in ${\rm dcl}_{\mathbb{M}_{|I_v}}([c\,c']\,\mathbb{P}_{|I_v})$. But $[c\, c']=\{(c,c')\}$  and $c, c', v \in {\rm dcl}([a]\,\mathbb{P})$. It follows that the element of $J$ in ${\rm dcl}_{\mathbb{M}_{|I_v}}([c\,c']\,\mathbb{P}_{|I_v})$ is in ${\rm dcl}([a]\,\mathbb{P})$. 

In the case $c'$ is one of the end points of $I_v$, by a similar argument,  we get an element of $J$ in ${\rm dcl}_{\mathbb{M}_{|I_v}}([c]\,\mathbb{P}_{|I_v})$ which will also be in ${\rm dcl}([a]\,\mathbb{P})$. So in both case there is an  element of $X$ in ${\rm dcl}([a]\,\mathbb{P})$.

\qed \\

\end{section}

\begin{section}{Concluding remarks}\label{section rmks}

Note that some authors might define definable Skolem functions/definable choice in a slightly different way, explained below. The way we defined these notions in Section \ref{section prelimp def skolem} corresponds to the way they were used in the papers  developing definable analogues of algebraic topology tools mentioned in the Introduction. 

One may (and often people choose to do that) define  definable Skolem functions/definable choice in an  ``unparametric'' way: $\mathbb{M}$ has {\it unparametric definable Skolem functions} if for every $A$-definable family $\{X_t\}_{t\in  T}$ of nonempty definable subsets of some $M^k$, there is a $A$-definable function $f:T\to M^k$ such that $f(t)\in X_t$ for all $t\in T$; $\mathbb{M}$ has {\it unparametric definable choice}, if  furthermore, we can choose $f:T\to M^k$ such that   $f(t)=f(t')$ whenever $X_t=X_{t'}$. 

With these definitions it is straightforward to obtain the  following analogues of Propositions \ref{prop def skolem dcl} and  \ref{prop def choice dcl}.  

\begin{fact}\label{fact unparam}
If $\mathbb{M}$ is sufficiently saturated, then: 
\begin{itemize}
\item[$\bullet$]
 ${\mathbb M}$ has unparametric definable Skolem functions if and only if every  definable subset $\emptyset \neq X\subseteq M$ defined with parameters $a\in M^l$ has an element in $\dcl(a)$.

\item[$\bullet$]
${\mathbb M}$ has unparametric definable choice if and only if every  definable subset $\emptyset \neq X\subseteq M$ defined with parameters $a\in M^l$ has an element in $\dcl([a])$.
\end{itemize}
\end{fact}

Replacing the use of Propositions \ref{prop def skolem dcl} and  \ref{prop def choice dcl} throughout the paper by the two unparametric versions above, we obtain the following unparametric version of our main theorem:

\begin{thm}\label{thm main 1}
Let $\mathbb{M}$ be an o-minimal structure. Then the following are equivalent:
\begin{enumerate}
\item[(1)]
$\mathbb{M}$ has unparametric definable Skolem functions. 
\item[(2)]
There is a finite collection $\mathcal{Y}$ of $\emptyset $-definable open subintervals of $M$ such that:
\begin{itemize}
\item[-]
$M\setminus \bigcup \mathcal{Y}$ is a finite set of trivial points each defined over $\emptyset $.
\item[-]
each $Y\in \mathcal{Y}$ is the union of a uniformly $\emptyset $-definable family of group-intervals, each with a fixed positive element, parametrized by the end points of the intervals.
\item[-]
each $Y\in \mathcal{Y}$ has  a fixed $\emptyset $-definable element.
\end{itemize}
\item[(3)] $\mathbb{M}$  has unparametric definable choice.
\end{enumerate}
\end{thm}

Finally note also that our main results characterize (unparametric) definable Skolem functions/definable choice in a rather explicit way. An obvious question is whether there is such a characterization for (parametric) elimination of imaginaries. We wonder if the technics used here and the examples from \cite[Proposition 3.2]{Pi86} and counter-example from \cite{Johnson} might help  attack this question. 

\end{section}

\bibliography{References}
\bibliographystyle{plain}

\end{document}